\nonstopmode

\documentclass[reqno,final]{amsart}
\usepackage{amssymb,amsmath,euscript}

\usepackage{graphicx}
\usepackage{epsfig}
\usepackage{verbatim}
\usepackage[all]{xy}
\newdir{ >}{{}*!/-10pt/@{>}}
\newdir^{ (}{{}*!/-5pt/@^{(}}
\newdir_{ (}{{}*!/-5pt/@_{(}}
\def\cgaps#1{}
\def\Cgaps#1{}
\allowdisplaybreaks

\def\undersetbrace#1\to#2{\underbrace{#2}_{#1}}											   
\def\oversetbrace#1\to#2{\overbrace{#2}^{#1}}
\def\AMSunderset#1\to#2{\underset{#1}{#2}}
\def\AMSoverset#1\to#2{\overset{#1}{#2}}

\allowdisplaybreaks

\def\undersetbrace#1\to#2{\underbrace{#2}_{#1}}
\def\oversetbrace#1\to#2{\overbrace{#2}^{#1}}
\def\AMSunderset#1\to#2{\underset{#1}{#2}}
\def\AMSoverset#1\to#2{\overset{#1}{#2}}

\newtheorem{theorem}{Theorem}[section]

\newtheorem*{prop*}{Proposition}

\newtheorem*{thm*}{Theorem}

\newtheorem*{lem*}{Lemma}

\newtheorem*{cor*}{Corollary}

\newtheorem*{conj*}{Conjecture}

\newtheorem{lemma}[theorem]{Lemma}

\newtheorem{corollary}[theorem]{Corollary}
\newtheorem{proposition}[theorem]{Proposition}

\def\x{\times}
 \def\ov{\overline}
\def\wt{\widetilde}
 \newcommand{\rk}{\mathop {\mathrm {rk}}\nolimits}
\renewcommand{\Im}{\mathop {\mathrm {Im}}\nolimits}

\begin{document}

\def\O{\mathrm{O}}
\def\ASp{\mathrm{ASp}}
\def\Sp{\mathrm{Sp}}
\def\U{\mathrm{U}}
\def\GL{\mathrm{GL}}
\def\SO{\mathrm{SO}}
\def\SL{\mathrm{SL}}

\def\kappa{\varkappa}
\def\phi{\phi}

\def\R{\mathbb{R}}
\def\C{\mathbb{C}}

\def\la{\langle}
\def\ra{\rangle}

 \def\cA{\mathcal A}
\def\cB{\mathcal B}

\def\cC{\mathcal C}

\def\cD{\mathfrak D}

\def\cE{\mathcal E}
\def\cF{\mathcal F}
\def\cG{\mathcal G}
\def\cH{\mathcal H}
\def\cJ{\mathcal J}
\def\cI{\mathcal I}
\def\cK{\mathcal K}
 \def\cL{\mathcal L}
\def\cM{\mathcal M}
\def\cN{\mathcal N}
 \def\cO{\mathcal O}
\def\cP{\mathcal P}
\def\cQ{\mathcal Q}
\def\cR{\mathcal R}
\def\cS{\mathcal S}
\def\cT{\mathcal T}
\def\cU{\mathcal U}
\def\cV{\mathcal V}
 \def\cW{\mathcal W}
\def\cX{\mathcal X}
 \def\cY{\mathcal Y}
 \def\cZ{\mathcal Z}
\def\0{{\ov 0}}
 \def\1{{^1}}
 \def\frA{\mathfrak A}
 \def\frB{\mathfrak B}
\def\frC{\mathfrak C}
\def\frD{\mathfrak D}
\def\frE{\mathfrak E}
\def\frF{\mathfrak F}
\def\frG{\mathfrak G}
\def\frH{\mathfrak H}
\def\frI{\mathfrak I}
 \def\frJ{\mathfrak J}
 \def\frK{\mathfrak K}
 \def\frL{\mathfrak L}
\def\frM{\mathfrak M}
 \def\frN{\mathfrak N} \def\frO{\mathfrak O} \def\frP{\mathfrak P} \def\frQ{\mathfrak Q} \def\frR{\mathfrak R}
 \def\frS{\mathfrak S} \def\frT{\mathfrak T} \def\frU{\mathfrak U} \def\frV{\mathfrak V} \def\frW{\mathfrak W}
 \def\frX{\mathfrak X} \def\frY{\mathfrak Y} \def\frZ{\mathfrak Z} \def\fra{\mathfrak a} \def\frb{\mathfrak b}
 \def\frc{\mathfrak c} \def\frd{\mathfrak d} \def\fre{\mathfrak e} \def\frf{\mathfrak f} \def\frg{\mathfrak g}
 \def\frh{\mathfrak h} \def\fri{\mathfrak i} \def\frj{\mathfrak j} \def\frk{\mathfrak k} \def\frl{\mathfrak l}
 \def\frm{\mathfrak m} \def\frn{\mathfrak n} \def\fro{\mathfrak o} \def\frp{\mathfrak p} \def\frq{\mathfrak q}
 \def\frr{\mathfrak r} \def\frs{\mathfrak s} \def\frt{\mathfrak t} \def\fru{\mathfrak u} \def\frv{\mathfrak v}
 \def\frw{\mathfrak w} \def\frx{\mathfrak x} \def\fry{\mathfrak y} \def\frz{\mathfrak z} \def\frsp{\mathfrak{sp}}
 \def\bfa{\mathbf a} \def\bfb{\mathbf b} \def\bfc{\mathbf c} \def\bfd{\mathbf d} \def\bfe{\mathbf e} \def\bff{\mathbf f}
 \def\bfg{\mathbf g} \def\bfh{\mathbf h} \def\bfi{\mathbf i} \def\bfj{\mathbf j} \def\bfk{\mathbf k} \def\bfl{\mathbf l}
 \def\bfm{\mathbf m} \def\bfn{\mathbf n} \def\bfo{\mathbf o} \def\bfp{\mathbf p} \def\bfq{\mathbf q} \def\bfr{\mathbf r}
 \def\bfs{\mathbf s} \def\bft{\mathbf t} \def\bfu{\mathbf u} \def\bfv{\mathbf v} \def\bfw{\mathbf w} \def\bfx{\mathbf x}
 \def\bfy{\mathbf y} \def\bfz{\mathbf z} \def\bfA{\mathbf A} \def\bfB{\mathbf B} \def\bfC{\mathbf C} \def\bfD{\mathbf D}
 \def\bfE{\mathbf E} \def\bfF{\mathbf F} \def\bfG{\mathbf G} \def\bfH{\mathbf H} \def\bfI{\mathbf I} \def\bfJ{\mathbf J}
 \def\bfK{\mathbf K} \def\bfL{\mathbf L} \def\bfM{\mathbf M} \def\bfN{\mathbf N} \def\bfO{\mathbf O} \def\bfP{\mathbf P}
 \def\bfQ{\mathbf Q} \def\bfR{\mathbf R} \def\bfS{\mathbf S} \def\bfT{\mathbf T} \def\bfU{\mathbf U} \def\bfV{\mathbf V}
 \def\bfW{\mathbf W} \def\bfX{\mathbf X} \def\bfY{\mathbf Y} \def\bfZ{\mathbf Z} \def\bfw{\mathbf w}
 \def\R {{\mathbb R }} \def\C {{\mathbb C }} \def\Z{{\mathbb Z}} \def\H{{\mathbb H}} \def\K{{\mathbb K}}
 \def\N{{\mathbb N}} \def\Q{{\mathbb Q}} \def\A{{\mathbb A}} \def\T{\mathbb T} \def\P{\mathbb P} \def\G{\mathbb G}
 \def\bbA{\mathbb A} \def\bbB{\mathbb B} \def\bbD{\mathbb D} \def\bbE{\mathbb E} \def\bbF{\mathbb F} \def\bbG{\mathbb G}
 \def\bbI{\mathbb I} \def\bbJ{\mathbb J} \def\bbL{\mathbb L} \def\bbM{\mathbb M} \def\bbN{\mathbb N} \def\bbO{\mathbb O}
 \def\bbP{\mathbb P} \def\bbQ{\mathbb Q} \def\bbS{\mathbb S} \def\bbT{\mathbb T} \def\bbU{\mathbb U} \def\bbV{\mathbb V}
 \def\bbW{\mathbb W} \def\bbX{\mathbb X} \def\bbY{\mathbb Y} \def\kappa{\varkappa} \def\epsilon{\varepsilon}
 \def\phi{\varphi} \def\le{\leqslant} \def\ge{\geqslant}

\def\O{\mathrm O}
\def\cAL{\mathcal{AL}}
\def\Symm{\mathrm {Symm}}
\def\USymm{\mathrm {USymm}}

\def\rS{\mathrm {S}}
\def\rB{\mathrm {B}}

\def\sm{\smallskip}
\def\wh{\widehat}
\def\tr{\mathrm{tr\,}}
\def\we{\mathrm{we}}
\def\on{\operatorname}

\title
{The Lagrangian Radon Transform and the Weil representation
}
\author{Giuseppe Marmo, Peter W. Michor, Yury A. Neretin}
\address{
G.\  Marmo: Dipartimento di Scienze Fisice, Universit\'a di Napoli
Federico II and INFN, Sezione di Napoli, via Cintia, 80126 Napoli, Italy }
\email{marmo@na.infn.it}
\address{
Peter W. Michor:
Fakult\"at f\"ur Mathematik, Universit\"at Wien,
Oskar-Morgenstern-Platz 1,, A-1090 Wien, Austria} 
\email{Peter.Michor@univie.ac.at}
\address{
Yury A. Neretin:
Fakult\"at f\"ur Mathematik, Universit\"at Wien,
Oskar-Morgenstern-Platz 1, A-1090 Wien, Austria. \&
Institute for Theoretical and Experimental Physics, Moscow, Russia
\& MechMath Department, Moscow State University
}
\email{
neretin@mccme.ru}
\date{\today}
\thanks{
YuAN was supported by FWF Projects P~22122 and P~25142}

\begin{abstract} 
We consider the operator $\cR$, which sends a function on $\R^{2n}$ to its integrals over all 
affine Lagrangian subspaces in $\R^{2n}$. We discuss properties of the operator $\cR$ and of
the representation of the affine symplectic group in several function spaces on $\R^{2n}$.
\end{abstract}
\def\LaTeXonly{}

\maketitle

\section{Introduction}

Quantum tomography \cite{MMS06, AFM12} has lately emerged as an alternative picture of
quantum mechanics \cite{IMM11, IMM12, MMS05}, alternative to the Dirac-Schr\"odinger picture on
Hilbert spaces or the $C^*$-algebraic description of Heisenberg and the
Wigner-Weyl formalism on phase-space.

The main aim of this paper is to provide a geometrical framework for the
transition from quasi-distributions, in the Weyl-Wigner description of
quantum mechanics on phase-space, to
fair probability distributions \cite{MMS08}, in the tomographic description of
quantum mechanics.

From the point of view of the present paper, the essential idea, behind
quantum tomography, may be illustrated by a simple example. We consider the
space of lines in the X-rays Radon transform as the space of Lagrangian
subspaces when the starting vector space is considered to be a two
dimensional symplectic vector space (phase-space). We may figure
out from here that a possible connection between the standard Radon
transform and the formulation of quantum mechanics on phase space might
exist. Adding the symplectic ingredient to the starting vector space makes
available a Fourier symplectic transform and allows to consider the family of
Lagrangian subspaces and, what is more relevant for the physical
interpretation, a Weyl system \cite{EMS04}. We recall that a Weyl system is a
projective unitary representation of the abelian vector group with the
symplectic structure playing the role of the cocycle in the projective
unitary representation. 

Roughly speaking, a state, considered as a rank-one projector, is associated,
by means of the Weyl operators \cite{MMM00}, with a function on phase-space, this
construction is often called "dequantization". Due to uncertainty relations
this function in generic cases will take both positive and negative values
and therefore cannot be considered to represent a fair probability
distribution. However, after integration on Lagrangian subspaces, the
resulting function will be a fair probability distribution (usually called
the marginal distribution) as a function on the space of Lagrangian
subspaces. The aim of quantum tomography is to reconstruct the original state
when a sufficient family of fair probability distributions is given. In
general these probability distributions may be replaced with a sufficient
set \cite{MMS06b} of expectation value functions or a quorum of operators. 

These expectation value functions are considered to be directly accessible
to an experimental determination and this makes more appealing the approach
also from the point of view of quantum field theory \cite{ILM12}. 

The reconstruction of the states amounts to the inversion of the Radon
transform first to obtain a function on phase space, thereafter we perform  a
further inversion from functions to operators by means of the Weyl map. A
relevant point is to identify those probability distribution functions on
the manifold of Lagrangian subspaces which are candidates to represent
quantum states once the inverse Radon transform is applied to them. 

The present paper will be mostly concerned with the definition of the Radon
transform on Lagrangian subspaces and its inversion, making it ready for use
both in classical and quantum tomography. In a coming paper we shall consider
the quantum framework more closely. 

\section{Introducing the Lagangian Radon transform }

\subsection*{The Lagrangian Radon transform}
Consider the space $\R^{2n}$ equipped with the standard skew-symmetric non-degenerate bilinear form,
the group $\ASp(2n,\R)$ of affine symplectic transformations, and the representation
of this group in functions on $\R^{2n}$.

The
Lagrangian Radon transform $\cR$ is defined in the following way.
Let $f$ be a suitable function on $\R^{2n}$. 
For any affine Lagrangian subspace $L\subset \R^{2n}$ we consider the integral
$$
\cR f(L)= \int_L f(x) dx.
$$
This integral transformation is similar to the usual Radon transform (i.e., integration over
arbitrary $k$-dimensional planes in $\R^m$),
see, e.g. \cite{Gra}, \cite{GGG},  \cite{Helgason}, \cite{Helgason94}.
Various operators of this
type  were considered in the literature, see \cite{Petrov67}, \cite{Petrov70},
\cite{Nievergelt86}, \cite{OurnychevaRubin08}, \cite{GonzalezKakehi06}.
The Lagrangian Radon transform was discussed by Grinberg \cite{Grinberg85}.
Nevertheless, we return to this old problem. 
Our main interest is the  relation of the Lagrangian Radon transform
with the Weil representation.

{\it We suppose $n>1$, otherwise some our statements become incorrect.}

\subsection*{The Determinantal system $\cD_{flat}$ of PDE} Decompose $\R^{2n}$ as a direct sum $\R^n\oplus\R^n$
with coordinates $(x,y)$. An affine Lagrangian subspace in general position has the form
$$x=Ty+\tau,$$
where $T$ is a symmetric $n\times n$ matrix 
and $\tau\in\R^n$. Consider the following matrix
\begin{equation}\label{determinantalPDE}
\begin{pmatrix}
2\frac\partial{\partial t_{11}} & \frac\partial{\partial t_{12}} & \dots&
     \frac\partial{\partial t_{1n}} &\frac\partial{\partial \tau_1}
\\
\frac\partial{\partial t_{12}} & 2\frac\partial{\partial t_{22}} & \dots&
   \frac\partial{\partial t_{2n}} &\frac\partial{\partial \tau_2}
\\
\vdots&\vdots& \ddots&\vdots&\vdots
\\
\frac\partial{\partial t_{1n}} & \frac\partial{\partial t_{2n}} & \dots &
          2 \frac\partial{\partial t_{nn}} &\frac\partial{\partial \tau_n}
\\
\frac\partial{\partial \tau_1}&\frac\partial{\partial \tau_2}&\dots& 
\frac\partial{\partial \tau_n}&0
\end{pmatrix}
\end{equation}
Any $C^1$-function of the form $\cR f$  is annihilated by all $3\times 3$ minors 
of this matrix (these equations were derived by Greenberg \cite{Grinberg85}).

Let $g$ satisfy this system of PDE.
If $g$ has  'good'
behavior at infinity (see (\ref{eq:boundary})),
 then it has form $g=\cR f$ (Theorem \ref{th:image},
see the corresponding statement for the usual Radon transform in \cite{Gra}, Theorem 5).
However, the determinantal PDE system \eqref{determinantalPDE} has 'extraneous' solutions
(see Proposition \ref{pr:extraneous} and Section 8). 

We also describe the image of the space of solutions of the determinantal system under the Fourier 
transform and we show the invariance of local solutions with respect to the natural
changes of variables (Theorem \ref{th:invariance}); this is one of numerous 
'higher' analogs of the conformal invariance of the Laplace equation (see \cite{Neretin11}, Section 7.8).

\subsection*{Relation with the Weil representation}
It is more-or-less evident, that  the representation of $\ASp(2n,\R)$ in $L^2(\R^{2n})$ is a tensor product of
the Weil representation of $\ASp(2n,\R)$ and its dual representation
(see Section 9). This allows to construct a continuation of functions of the form $\cR f$ to a certain complex domain,
see Section 10.

On the other hand, there are well known realizations of degenerate highest weight representations
of semisimple Lie groups
in spaces of distributions supported by determinantal varieties; see
Vergne, Rossi \cite{VR} and see also an exposition in \cite[Section 7.5]{Neretin11}. We imitate 
this construction
for the representation of $\ASp(2n,\R)$ in $L^2(\R^{2n})$; see Section 8. We get a certain space of distributions
supported by a real algebraic variety, that contains two  mustaches (similar to the Whitney umbrella).

\subsection*{The group $\ASp(2n,\R)$}
The group $\ASp(2n,\R)$ is a member of a class of non-semisimple groups which are close to 
semisimple ones.
Such groups (and the group $\ASp(2n,\R)$ itself) were the topic of several works; see
\cite{BerndtSchmidt98}, \cite{Berceanu07}, \cite{KOP11}. Some standard  statements
about representations of semi-simple groups become non-obvious or even false for such groups
(for instance, the Radon transform is an intertwining operator, which is defined on Schwartz space and
does not admit a continuation to the Hilbert space $L^2(\R^{2n})$).
On the other hand some new elements appear, see an unusual 
expression for an intertwining operator in Theorem \ref{th:intertwining}.

\subsection*{Notation}
Let $A$ be a matrix. 
\begin{itemize}
	\item 
$A^t$ is the transposed matrix.
\item
$\ov A$ is the matrix obtained from $A$ by coefficient-wise complex  conjugation.
\item
$A^*=\ov A^t$ is the adjoint matrix.
\item
$\bbI=\bbI_n$ is the unit matrix of order $n$.
\item
$\Symm(n)$ is the set of symmetric matrices of order $n$.
\end{itemize}

\section{The symplectic group}

\subsection*{The complex symplectic group}
Consider the linear  space $W=\C^{2n}=\C^n\oplus \C^n$. Equip it with
the skew-symmetric bilinear form $\Lambda$ defined by the matrix
\begin{equation}
\begin{pmatrix}
 0&\bbI\\ -\bbI&0
\end{pmatrix}
\end{equation}
We denote by $\Sp(2n,\C)$ the {\it complex symplectic group}, i.e., the group of all matrices preserving the form
$\Lambda$. We have:
\begin{equation}
g=\begin{pmatrix}
    a&b\\c&d
   \end{pmatrix}\in \Sp(2n,\C) \iff
 \begin{pmatrix}
    a&b\\c&d
   \end{pmatrix}
\begin{pmatrix}
 0&\bbI\\ -\bbI&0
\end{pmatrix}
\begin{pmatrix}
    a&b\\c&d
   \end{pmatrix}^t
=
\begin{pmatrix}
 0&\bbI\\ -\bbI&0
\end{pmatrix}
\label{eq:spdef}
\end{equation}
We also consider the affine symplectic group $\ASp(2n,\C)$ 
of affine symplectic transformations of $\C^{2n}$. It is the semidirect product of
$\Sp(2n,\C)$ with the translation group as normal subgroup, 
$$\ASp(2n,\C)=\Sp(2n,\C)\ltimes \C^{2n}.$$
It is convenient to realize $\ASp(2n,\C)$ as the group of matrices
$$
g=\begin{pmatrix}
    a&b&r\\c&d&s\\0&0&1
   \end{pmatrix}
\text{   of size }n+n+1, \text{  where  }
\begin{pmatrix}
    a&b\\c&d
   \end{pmatrix}\in \Sp(2n,\C).
$$ 
Its action on $\C^{2n}$ is given by the formula
$$
\begin{pmatrix}
    a&b&r\\c&d&s\\0&0&1
   \end{pmatrix}
\begin{pmatrix}
x\\ y\\1
\end{pmatrix}
=
\begin{pmatrix}
    a x +  by +r\\cx +  dy  +s\\1
   \end{pmatrix}
.$$

\subsection*{The real symplectic group. Two models}
 For details, see \cite{Neretin11}, Section 3.2.
The real symplectic group $\Sp(2n,\R)$ is the group of real 
matrices satisfying the same condition (\ref{eq:spdef}).

Sometimes it is convenient to realize $\Sp(2n,\R)$ as the group
of matrices of the form 
\begin{equation}
g=\begin{pmatrix}
 \Phi&\Psi\\
\ov \Psi&\ov\Phi
\end{pmatrix}
,
\label{eq:PhiPsi}
\end{equation}
which are contained in the complex symmetric group $\Sp(2n,\C)$, i.e., they satisfy the condition
\begin{equation}
g\begin{pmatrix} 0&\bbI\\-\bbI&0\end{pmatrix}g^t=
\begin{pmatrix} 0&\bbI\\-\bbI&0\end{pmatrix}
\label{eq:PhiPsiSp}
.
\end{equation}
For matrices of the form (\ref{eq:PhiPsi}) we can replace the  condition (\ref{eq:PhiPsiSp}) by the equivalent condition
\begin{equation}
g\begin{pmatrix} \bbI & 0\\0 & -\bbI\end{pmatrix}g^*=
\begin{pmatrix} \bbI & 0\\0 &-\bbI\end{pmatrix}
\label{eq:Unn}
.
\end{equation}

\noindent
{\sc Remarks.} a) Matrices of the form
(\ref{eq:PhiPsi}) preserve the {\it real} linear subspace in $\C^{2n}$ consisting 
of vectors $\begin{pmatrix}h\\ \ov h \end{pmatrix}$, where $h$ ranges in $\C^n$.

\sm

b) Condition (\ref{eq:Unn}) means that matrices of the form
(\ref{eq:PhiPsi}) preserve the 
Hermitian form with matrix $\begin{pmatrix} \bbI & 0\\0 & -\bbI\end{pmatrix}$.
Also, the pair of equations (\ref{eq:PhiPsiSp})--(\ref{eq:Unn})
implies that $g$ has block structure 
$\begin{pmatrix}
 \Phi&\Psi\\
\ov \Psi&\ov\Phi
\end{pmatrix}$.
\hfill $\square$

\sm

The passage from one model to another is given by conjugation
 \begin{equation}
g\mapsto K g K^{-1},
\qquad\text{where $K=\frac 1{\sqrt 2}
\begin{pmatrix}
\bbI&i\cdot\bbI\\i\cdot\bbI&\bbI 
\end{pmatrix}$}
\label{eq:link}
\end{equation}
Notice that $K\in\Sp(2n,\C)$.

\sm

Consider the affine real symplectic group
$$
\ASp(2n,\R)=\Sp(2n,\R)\ltimes \R^{2n}
.$$
We realize it as a subgroup in $\ASp(2n,\C)$ in two ways:
\begin{itemize}
	\item 
As the group of real matrices of the form
\begin{equation}
\begin{pmatrix}
    a&b&r\\c&d&s\\0&0&1
   \end{pmatrix}, \qquad \text{where $\begin{pmatrix}
    a&b\\c&d
   \end{pmatrix}\in \Sp(2n,\R)$ } 
   \label{eq:aff-r}
   .
   \end{equation}
\item
   Or the group of complex matrices of the form
  \begin{equation}
  \begin{pmatrix}
    \Phi&\Psi&h\\\ov\Psi&\ov \Phi&\ov h\\0&0&1
 \end{pmatrix},\qquad \text{where $\begin{pmatrix}
    \Phi&\Psi\\\ov\Psi&\ov \Phi
 \end{pmatrix}\in\Sp(2n,\C)$}  
 .
 \label{eq:aff-c}
 \end{equation}
\end{itemize}
 
 \subsection*{Subgroups of $\Sp(2n,\R)$}
 We consider the following subgroups in $\Sp(2n,\R)$:
\begin{itemize}
  \item 
The subgroup $\GL(n,\R)$ consisting of matrices
of the form $\begin{pmatrix}a&0\\0&(a^t)^{-1} \end{pmatrix}$, where $a$ is an arbitrary
invertible matrix.
 
\item
The abelian subgroup $N_+$ consisting of matrices
 $\begin{pmatrix}\bbI&b\\0&\bbI \end{pmatrix}$.%

\item
The abelian subgroup $N_-$ consisting of matrices
 $\begin{pmatrix}\bbI&0\\c&\bbI \end{pmatrix}$, where $c$ is a symmetric matrix.
\end{itemize}

 The product $N_-\GL(n,\R) N_+$ is a dense subset in $\Sp(2n,\R)$;
 see, e.g., \cite[formula (1.2.14)]{Neretin11}.

 Since
 $$
 \begin{pmatrix} 0&\bbI\\-\bbI&0\end{pmatrix}
  \begin{pmatrix} \bbI&b\\0&\bbI\end{pmatrix}
  \begin{pmatrix} 0&\bbI\\-\bbI&0\end{pmatrix}^{-1}=
    \begin{pmatrix} \bbI&0\\-b&\bbI\end{pmatrix}
 ,$$
 the subgroups $N_+$, $\GL(n,\R)$, and the element  $\begin{pmatrix} 0&\bbI\\-\bbI&0\end{pmatrix}$
 together generate the group $\Sp(2n,\R)$.

 We also define a {\it parabolic subgroup} $P=\GL(n,\R)\ltimes N_+$,
 it consists of symplectic matrices of the form $\begin{pmatrix}a&b\\0&(a^t)^{-1} \end{pmatrix}$.
 
 Next,  $\Sp(2n,\R)$ contains the unitary group $\U(n)$. Indeed, consider the Euclidean space $\C^n$
 with the standard inner product $\la\cdot,\cdot\ra$. Let us regard $\C^{n}$ as a real
 linear space $\R^{2n}$. Then $\Im\la\cdot,\cdot\ra$ is a skew symmetric 
 bilinear form and unitary operators in $\C^n$ preserve it.
 
 In the first model the subgroup $\U(n)\subset \Sp(2n,\R)$ is
 given by symplectic matrices of the form
 $\begin{pmatrix}a&b\\-b&a\end{pmatrix}$.
 In the second model, by
 $\begin{pmatrix}\Phi&0\\0&\ov\Phi \end{pmatrix}$, where $\Phi$ is unitary. 
 
 \sm
 
 Finally, we define the subgroup $\O(n)=\GL(n,\R)\cap \U(n)$. In both models
 it is given by matrices of the form  $\begin{pmatrix}a&0\\0&a\end{pmatrix}$,
 where $a$ is a real orthogonal matrix.

 \subsection*{Lagrangian Grassmannian%
 \label{ss:grass}}
 For details, see \cite{Neretin11}, Section 3.3.
 Recall that we equipped $\C^{2n}$ with a non-degenerate skew-symmetric bilinear form 
 $\Lambda$.  A subspace $L$  is called {\it isotropic} if $\Lambda=0$ on $L$.
 A subspace $L$ is called {\it Lagrangian} if $L$ is an isotropic subspace of maximal
 possible dimension (i.e., $\dim L=n$). Denote by $\cL(n,\C)$ the set of all 
 Lagrangian subspaces in $\C^{2n}$. 
 
 Consider the following subspaces
 in $W:=\C^{2n}=\C^n\oplus \C^n$:
 $$
 V_+=\C^n\oplus 0,\qquad V_-=0\oplus \C^n
 $$
Let $T$ be a symmetric ($T=T^t$) operator 
$V_-\to V_+$. Then its graph $graph(T)$ is a Lagrangian
subspace in $V_-\oplus V_+= \C^{2n}$. Conversely,
any Lagrangian subspace $L\subset\C^{2n}$ such that $L\cap V_+=0$
has such form. Thus we identify $\Symm(n,\C)$
with an open dense subset in $\cL(n,\C)$.

The action of $\Sp(2n,\C)$ on $\cL(n,\C)$ in terms of these coordinates is given by
\begin{equation}
\begin{pmatrix}a&b\\c&d\end{pmatrix}:
\quad T\mapsto (aT+b)(cT+d)^{-1}
\label{eq:lin-frac}
.
\end{equation}
Indeed, consider the linear subspace $L$ consisting of vectors $(Ty,y)\in \C^n\oplus \C^n$.
Then $gL$ consists of vectors
\begin{equation}
\begin{pmatrix}
a&b\\c&d
\end{pmatrix}
\begin{pmatrix}
Ty\\y
\end{pmatrix}=
\begin{pmatrix}(aT+b)y\\ (cT+d)y \end{pmatrix}
\label{eq:trv}
\end{equation}
and we pass to the new variable $z=(cT+d)y$.
\hfill $\square$

\sm

Note that this construction is not canonical:
we fixed the decomposition $W=V_-\oplus V_+$. However, we can replace 
$V_-$ and $V_+$ by an arbitrary pair of transversal Lagrangian subspaces $Y$, $Z$. 
The bilinear form $\Lambda$ determines a nondegenerate pairing
$Y\times Z\to \C$ by $(y,z)\mapsto \Lambda(y,z)$, therefore $Z$
is canonically isomorphic to the dual space $Y'$ of $Y$. Hence
for any operator $T:Y\to Z$ the transposed operator $T^t$ is also an operator
$Y\to Z$. Therefore the condition $T=T^t$ makes sense. As above, for any symmetric operator 
$T:Y\to Z$ we get a Lagrangian subspace. 
This describes an atlas on $\cL(n,\C)$; charts are enumerated by pairs of transversal Lagrangian subspaces.

 Our topic is the real Lagrangian Grassmannian $\cL(n,\R)\subset \cL(n,\C)$.
 We can define coordinates on $\cL(n,\R)$ in two ways.
 
First, we can repeat literally the construction given above.
Thus we identify $\Symm(n,\R)$
with an open dense subset in $\cL(n,\R)$. And considering all pairs of transversal Lagrangian 
subspaces instead of $V_-, V_+$, we get an atlas of $\cL(n,\R)$. 
 
 Since
 $V_-$ is invariant with respect to the parabolic subgroup
 $P$, we have
 $$
 \cL(n,\R)\simeq\Sp(2n,\R)/P.
 $$
 
 However, there is an another way. We can complexify the real symplectic
 space $\R^{2n}$ and choose Lagrangian subspaces in the complexification.
 The simplest way is to take two transversal Lagrangian subspaces corresponding to $T=i\bbI, 
 -i\bbI$, namely $W_-=graph (i\bbI)$, $W_+=graph(-i\bbI)$.
 Notice, that complexifying a real skew-symmetric bilinear form $\Lambda$
 we get {\it two forms}, the first is the bilinear extension of $\Lambda$
 (we preserve the same notation $\Lambda$):
 \begin{equation}
 \Lambda(v_1+iv_2,w_1+iw_2)=
 \Lambda(v_1,w_1)-\Lambda(v_2,w_2)+i\Lambda(v_1,w_2)+i\Lambda(v_2,w_1)
 \label{eq:form-Lambda}
 \end{equation}
 The second is the sesquilinear extension
 \begin{equation}
 \label{eq:M}
 M(v,w)=\frac 1i \Lambda(v,\ov w)
.
\end{equation}
 We divide the form by $i$ to make it Hermitian,
 $M(v,w)=\ov{M(w,v)}$. The form  $M$ is positive definite on $W_-$ and negative definite
 on $W_+$. On the other hand, a complexification of a real subspace must be isotropic with
 respect to both forms, $\Lambda$ and  $M$.  This easily  implies that in our case $T$ is unitary.
 
 Thus we identified $\cL(n,\R)$ with the set of all symmetric unitary matrices.
 Notice that $\cL(n,\R)$ is a $\U(n)$-homogeneous space (see also \cite{Michor08}, Section 31.7), 
 $$\cL(n,\R)\simeq\U(n)/\O(n).$$
 
 The two systems of coordinates are related by the Cayley transform
 \begin{equation}
 T=(S+i)(iS+1)^{-1}, \quad\text{ where }T=T^t\text{  is real and }S=S^t\text{  is unitary.}
 \label{eq:TS}
 \end{equation}

 \subsection*{Affine Lagrangian Grassmannian}
 Consider the set $\cAL(n,\C)$ of all affine subspaces in $\C^{2n}$ obtained by translates
 of Lagrangian linear subspace. A  subspace in general position can be given by the equation
\begin{equation}
 x=Ty+\tau, \quad \text{ where }(x,y)\in \C^{2n},\quad T\in\Symm(n,\C),\quad \tau\in\C^n.
 \label{eq:aff-T-tau}
 \end{equation}
 
 \begin{lemma}
 The action of $\ASp(2n,\R)$ on $\cAL(n,\C)$ on the dense open subset described above is given by the formula
 \begin{multline}
 \begin{pmatrix}
    a&b&r\\c&d&s\\0&0&1
   \end{pmatrix}:\quad (T,\tau)
   \mapsto \\ \mapsto \bigl((aT+b)(cT+d)^{-1},\quad  a\tau+r\, -(aT+b)(cT+d)^{-1} (c\tau+s)  \bigr)
\label{eq:action}
.
   \end{multline}
 \end{lemma}
 
 {\sc Proof.}
We write
$$
\begin{pmatrix}a&b&r\\
 c&d&s\\
0&0&1
\end{pmatrix}
\begin{pmatrix}
 Ty + \tau\\ y\\1
\end{pmatrix}
=
\begin{pmatrix}
 (aT+b)y +a\tau+r
\\
(cT+d)y+c\tau+s
\\
1
\end{pmatrix}
$$
and pass to the variable $z= (cT+d)y+c\tau+s$.
\hfill $\square$
 
 \sm
 
 The expression in the right hand side of (\ref{eq:action}) admits a non-obvious transformation:
 
 \begin{lemma}
 If $\begin{pmatrix}a&b\\c&d \end{pmatrix}\in\Sp(2n,\C)$ and $T\in \Symm(n,\C)$, then
 \begin{equation}
 a-(aT+b)(cT+d)^{-1} c= (Tc^t+d^t)^{-1}=\big((cT+d)^t\big)^{-1}
 \label{eq:t}
 \end{equation}
 whenever one side makes sense.
 \end{lemma}
 
{\sc Proof.} For generic elements we compute as follows. 
\begin{multline*}
a-(aT+b)(cT+d)^{-1} c
=\\= \bigl(ac^{-1}(cT+d)-(aT+b)\bigr)(cT+d)^{-1} c
=(ac^{-1}d-b)(cT+d)^{-1} c
.
\end{multline*}
Next, by (\ref{eq:spdef}), we have
$$
-dc^t+cd^t=0,\qquad ad^t-bc^t=1
.
$$
Therefore
$$
ac^{-1}d-b= a(c^{-1}d)-b=a d^t (c^t)^{-1} -b= (1+bc^t)(c^t)^{-1}-b=(c^t)^{-1}. 
$$
Finally,
$$
(c^t)^{-1} (cT+d)^{-1}c=(Tc^t+c^{-1}dc^t)^{-1}=(Tc^t+d^t)^{-1}
$$
Our calculation is valid for invertible $c$. However both sides
of (\ref{eq:t}) are
continuous on the common domain of definity. Therefore the equation
holds for all $c$.
\hfill $\square$

\sm

We return to the  group $\ASp(2n,\R)$.
For the  real realization of $\ASp(2n,\R)$ we have precisely the same formulas.
For the complex realization of $\ASp(2n,\R)$, we have the following:

\begin{lemma}
An affine Lagrangian subspace descibed as  $\{(w,v)\in W_-\times W_+: v=Sw+\sigma\}$, where
$S$ is an unitary symmetric matrix is
a complexification of a real affine Lagrangian subspace if and only if
$$
\sigma=-S\ov\sigma.
$$
\end{lemma}

{\sc Proof.} A subspace $L$ must be closed with respect to complex conjugation.
In our coordinates it is given by the formula
$(v,w)\mapsto(\ov w,\ov v)$. Therefore, for any $w$ there is $h$ such that
$$(Sw+\sigma,w)=(\ov h, \ov{Sh+\sigma}).$$
Therefore,
$$
S(\ov{Sh+\sigma})+\sigma=\ov h
$$
and this implies the required statement.
\hfill $\square$

\section{Radon transform. Determinantal system}

Denote by $\cS(\R^{2n})$ the Schwartz space on $\R^{2n}$ of rapidly decreasing functions.

\subsection*{Radon transform}
Consider the space $\R^{2n}$ equipped with the standard symplectic form $\Lambda$.
Consider the space $\wt \cAL(n,\R)$ whose points are pairs
\begin{multline*}
\text{(Lagrangian affine subspace $L$,}
\\ \text{a translation-invariant positive volume
 form $d\omega$ on $L$)}
\end{multline*}
The space $\wt\cAL(n,\R)$ is a fiber bundle over $\cAL(n,\R)$ with fiber $\R^+$.

For a rapidly decreasing smooth function $f\in \cS(\R^{2n})$ we define its {\it Radon transform} 
as the function
$$
\cR f(L,\omega)=\int_L f(v)\,d\omega(v).
$$
Obviously, $\cR f(L,\omega)$ is a smooth function on $\wt \cAL(n,\R)$ satisfying
\begin{equation}
\cR f(L,t\cdot\omega)=t\cdot\cR f(L,\omega),\qquad t>0
\label{eq:homogen}
\end{equation}

 \subsection*{Explicit formula. Variant 1.}
 Now we wish to specify a form $\omega$ for a fixed $L\in\cAL(n,\R)$.
 Introduce coordinates of the first kind on the Grassmannian; i.e.,
 we identify the set $\Symm(n,\R)\times\R^n$ with an open dense subset in  
 $\cAL(n,\R)$. Recall that an affine subspace $L$ is now the graph of
 an operator, $x=Ty+\tau$, where $y$ ranges in $\R^n$. Thus we have 
 the map $\R^n\to L$
 given by
 $$y\mapsto (Ty+\tau,y).$$
 We equip $L$ with the pushforward 
 $d\lambda_L$ of the Lebesgue measure on $\R^n$.
 
For $f\in\cS(\R^{2n})$ we define a function on $\Symm(n,\R)\times \R^n$
by
\begin{equation}
\cR_{flat} f(T,\tau)=\int_L f(v)\,d\lambda_L(v)=  \int_{\R^n} f(Ty+\tau,y)\,dy
\label{eq:R-flat}
\end{equation}

\begin{proposition}
The Radon transform $\cR_{flat}$ sends the action of $\ASp(2n,\R)$ on $\cS(\R^{2n})$
to the action $\rho_{flat}$ on the space of functions on $\cAL(n,\R)$, which on generic elements is 
given by:
 \begin{align}
  \label{eq:perenos-1}
&\rho_{flat}
\begin{pmatrix}a&b&r\\
 c&d&s\\
0&0&1
\end{pmatrix}
h(T,\tau)
=\\&=
h\bigl((aT+b)(cT+d)^{-1},  a\tau+r\, -(aT+b)(cT+d)^{-1} (c\tau+s)  \bigr)
 \bigl|\det (cT+d)\bigr|^{-1}
  .\notag
\end{align}
The formula holds for invertible $cT+d$. For $h=\cR_{flat}f$, $f\in\cS(\R^{2n})$,
the right-hand side admits a smooth continuation to the hypersurface
$\det(cT+d)=0$.
\end{proposition}
 
{\sc Proof.}
For $h=\cR_{flat}f$ and the matrix $g$ in the proposition, we must evaluate
$$
\cR_{flat} f(gL)=
\int_{\R^n} f\bigl( a(Ty+\tau)+by+r,c(Ty+\tau)+dy+s\bigr)\,dy
$$
 Substituting $z=(cT+d)y+s$, $y=(cT+d)^{-1}(z-s)$ we get (\ref{eq:perenos-1}). 
\hfill$\square$

\sm

{\sc Remark.} This formula shows that functions in the image of $\cR_{flat}$
satisfy rigid conditions at infinity. For instance, for $h(T,\tau)=\cR_{flat} f$ 
the function
\begin{equation}
h(T^{-1},-T^{-1}\tau) |\det(T)|^{-1}
\qquad
\text{is smooth.}
\label{eq:boundary}
\end{equation}

\subsection*{Explicit formula. Variant 2.}
Now we equip $\R^{2n}$ with the standard inner product.
Then any Lagrangian subspace $L$ is equipped with
a canonical surface measure $d\mu_L$ induced by the Euclidean structure.

\sm

{\sc Remark.} In particular,  each fiber of the bundle
$\wt \cAL\to \cAL$ has a canonical representative. Therefore the bundle is trivial.
\hfill $\square$

\begin{lemma}
Let $L$ be the graph of the operator $x=Ty+\tau$. Then
\begin{equation}
d\mu_L= \det(1+T^2)^{-1/2}d\lambda_L= |\det(1+iT)|^{-1}d\lambda_L
\label{eq:lambda-mu-1}
.
\end{equation}
If the complex coordinates of $L$ are given by $(S,\sigma)$, then
\begin{equation}
d\mu_L=\frac{\det|iS+1|}{2^n} d\lambda_L
.
\label{eq:lambda-mu-2}
\end{equation}

\end{lemma}

{\sc Proof.} Without loss of generality we can set $\tau=0$.
Consider the element of the subgroup $\U(n)\subset \Sp(2n,\R)$
given by the formula
$$
g=\begin{pmatrix} (1+T^2)^{-1/2}& T(1+T^2)^{-1/2}\\
-T(1+T^2)^{-1/2}&(1+T^2)^{-1/2}
 \end{pmatrix}
.$$
It sends the subspace $0\oplus\R^n$ to $L$,
$$
y\mapsto \bigl(T(1+T^2)^{-1/2}y, (1+T^2)^{-1/2}y\bigr)
$$
and we observe that measures $d\lambda_L$ and $d\mu_L$ differ by the factor
$\det (1+T^2)^{-1/2}$.

Next, $T$ and $S$ are related by
(\ref{eq:TS}),
and (\ref{eq:lambda-mu-1}) implies  (\ref{eq:lambda-mu-2}).
\hfill $\square$

\sm

Now we can define the following version of the Radon transform
$$
\cR_{comp} f(S,\sigma)=\int_L f(v)\,d\mu_L(v)
,
$$
where $(S,\sigma)$ are the complex coordinates of $L$ and where $f\in \cS(\R^{2n})$.

\begin{proposition}
The Radon transform $\cR_{comp}$ sends the action of $\ASp(2n,\R)$
 on $\cS(\R^{2n})$
to the action
\begin{multline}
\rho_{comp}
\begin{pmatrix}\Phi&\Psi&p\\
 \ov\Psi&\ov\Phi&\ov p\\
0&0&1
\end{pmatrix}
f(S,\sigma)=
\\
=
f\Bigl((\Phi S+\Psi)(\ov \Psi S +\ov\Phi)^{-1},  \Phi\sigma+p\, -(\Phi S+\Psi)(\ov\Psi S+\ov\Phi)^{-1} (\ov\Psi\sigma+\ov p)  \Bigr)
\cdot
\\ \cdot
 \bigl|\det (\ov \Psi S+\ov\Phi)\bigr|^{-1}
\label{eq:perenos-2}
\end{multline}
\end{proposition}

{\sc Proof.} In the big brackets we have a change of coordinates.
The expression for the last factor is a corollary of the following lemma.

\begin{lemma}
For $g=\begin{pmatrix}\Phi&\Psi\\\ov\Phi&\ov\Psi \end{pmatrix}$ we have 
$$
g_*d\mu_L=
|\det (\ov \Psi S+\ov\Phi)|^{-1}
d\mu_{gL}
.
$$
\end{lemma}

{\sc Proof.}
For a  matrix $g$ as given 
and a unitary symmetric matrix $S$
we define two expressions:
$$
c(g,S):=|\det(\ov\Psi S+\Phi)|,\qquad \wt c(g,S):=\frac{g_*d\mu_L}{d\mu_{gL}}.
$$
It can be checked that $c(g,T)$ satisfies the chain identity
\begin{equation}
c(g_1 g_2,T)=c\bigl(g_1, (\Phi_2 T+\Psi_2)(\ov\Psi_2 T+\Phi_2)^{-1}\bigr) c(g_2,T).
\label{eq:cocycle}
\end{equation}
By definition $\wt c(g,S)$ satisfies the same identity. Obviously, for
$h\in\U(n)$ we have
$$
c(h,S)=1\qquad \wt c(h,S)=1.
$$
Next we consider $p, q$ in the stabilizer $P$ of $V_-$ (the corresponding $S$ is 1); recall
that $P$ consists of matrices $\begin{pmatrix}a&b\\0&(a^t)^{-1} \end{pmatrix}$. By 
(\ref{eq:cocycle}), 
$$c(pq,1)=c(p,1)c(q,1), \qquad \wt c(pq,1)=\wt c(p,1)\wt c(q,1).
$$
All homomorphisms 
$P\to\R^\times_+$ equals $|\det (a)|^s$. 
It is sufficient to verify the identity $c(p,1)=\wt c(p,1)$
for $p=\begin{pmatrix}s\cdot \bbI&0\\0&s^{-1}\cdot\bbI \end{pmatrix}$;
this verification is straightforward.

Now we represent $g$ as $g=hp$ for $h\in\U(n)$ and $p\in P$; this is possible since the action of $\U(n)$
on $\cL$ is transitive. By  identity (\ref{eq:cocycle}), we have 
$$
c(hp,1)=c(p,1), \qquad \wt c(hp,1)=\wt c(p,1)
.
$$
This implies the statement.
\hfill $\square$

\sm

The next proposition establishes a link between the functions $\cR_{comp} f$ and $\cR_{flat} f$.

\begin{proposition}

\begin{multline}
\cR_{comp}f(S,\sigma)=
\\
=\cR_{flat}f\bigl((S+i)(iS+1)^{-1}, \sigma-(S+i)(iS+1)^{-1}i\sigma\bigr) \frac{2^n}{|\det(1+iS)|}
\label{eq:perenos-3}
\end{multline}

\end{proposition}

{\sc Proof.} We write the transformation (\ref{eq:action})
for the matrix 
$\begin{pmatrix}\bbI&i\cdot\bbI&0\\i\cdot\bbI&\bbI&0\\0&0&1 \end{pmatrix}$.
The factor $|\det(1+iS)|^{-1}$ arises from (\ref{eq:lambda-mu-2}).
\hfill $\square$

\subsection*{The system $\cD_{flat}$ of partial differential equations}
Consider the matrix composed of partial derivatives
\begin{equation}
\begin{pmatrix}
2\frac\partial{\partial t_{11}} & \frac\partial{\partial t_{12}} & \dots&
     \frac\partial{\partial t_{1n}} &\frac\partial{\partial \tau_1}
\\
\frac\partial{\partial t_{12}} & 2\frac\partial{\partial t_{22}} & \dots&
   \frac\partial{\partial t_{2n}} &\frac\partial{\partial \tau_2}
\\
\vdots&\vdots& \ddots&\vdots&\vdots
\\
\frac\partial{\partial t_{1n}} & \frac\partial{\partial t_{2n}} & \dots &
          2 \frac\partial{\partial t_{nn}} &\frac\partial{\partial \tau_n}
\\
\frac\partial{\partial \tau_1}&\frac\partial{\partial \tau_2}&\dots& 
\frac\partial{\partial \tau_n}&0
\end{pmatrix}
\label{eq:system}
\end{equation}

Each minor of this matrix determines a partial differential operator
 with constant coefficients.

\begin{theorem}
\label{th:system}
{\rm(}See {\rm \cite{Grinberg85}}{\rm)}
 For  smooth $f\in\cS(\R^{2n})$ the function 
$(T,\tau)\mapsto \cR_{flat} f(T,\tau)$ satisfies the system
of equations
\begin{equation}
D \,\cR_{flat} f(T)=0
\label{eq:system-2}
,
\end{equation}
where $D$ ranges in all $3\times 3$ minors of the matrix
{\rm(\ref{eq:system})}.
\end{theorem}

We denote this system of PDE by $\cD_{flat}$.

{\sc Proof.}
 We have 
\begin{align*}
\frac{\partial}{\partial t_{kl}}
\int_{\R^n} f(Ty+\tau,y)\,dy&=
\int_{\R^n} \left(y_k \frac\partial {\partial x_l}+ 
y_l \frac\partial {\partial x_k}  \right)f (Ty+\tau,y)\,dy,
\qquad \text{if $k\ne l$}
\\
\frac\partial{\partial t_{kk}}
\int_{\R^n} f(Ty+\tau,y)\,dy&=
\int_{\R^n} \left(y_k \frac\partial {\partial x_k} f \right) (Ty+\tau,y)\,dy,
\\
\frac\partial {\partial\tau_k}
\int_{\R^n} f(Ty+\tau,y)\,dy&=
\int_{\R^n} \left( \frac\partial {\partial x_k} f \right) (Ty+\tau,y)\,dy.
\end{align*}
Thus the matrix {\rm(\ref{eq:system})} acts as the following matrix of commuting differential operators on the functions 
$f$ inside the integral:
\begin{align*}
&\begin{pmatrix}
2x_1\frac{\partial}{\partial y_1} &(x_1\frac{\partial}{\partial y_2}+ x_2\frac{\partial}{\partial y_1})& \dots & (x_1\frac{\partial}{\partial y_n}+x_n\frac{\partial}{\partial y_1})& \frac{\partial}{\partial y_1}
\\
(x_1\frac{\partial}{\partial y_2}+ x_2\frac{\partial}{\partial y_1})& 2x_2\frac{\partial}{\partial y_2}& \dots & (x_2\frac{\partial}{\partial y_n}+x_n\frac{\partial}{\partial y_2})& \frac{\partial}{\partial y_2}
\\
\vdots&\vdots &\ddots &\vdots &\vdots 
\\
(x_1\frac{\partial}{\partial y_n}+ x_n\frac{\partial}{\partial y_1})& (x_2\frac{\partial}{\partial y_n}+x_n\frac{\partial}{\partial y_2})& \dots & 2x_n\frac{\partial}{\partial y_n}&\frac{\partial}{\partial y_n}
\\
\frac{\partial}{\partial y_1}&\frac{\partial}{\partial y_2}&\dots &\frac{\partial}{\partial y_n}  & 0
\end{pmatrix}=
\\&
=\begin{pmatrix}
x_1 \begin{pmatrix} \frac{\partial}{\partial y_1} \\ \vdots \\ \frac{\partial}{\partial y_n} \\ 0 \end{pmatrix}
+ \frac{\partial}{\partial y_1} \begin{pmatrix} x_1 \\ \vdots \\ x_n \\ 1 \end{pmatrix}, 
& \dots, &
x_n \begin{pmatrix} \frac{\partial}{\partial y_1} \\ \vdots \\ \frac{\partial}{\partial y_n} \\ 0 \end{pmatrix}
+ \frac{\partial}{\partial y_n} \begin{pmatrix} x_1 \\ \vdots \\ x_n \\ 1 \end{pmatrix}, 
& 1\begin{pmatrix} \frac{\partial}{\partial y_1} \\ \vdots \\ \frac{\partial}{\partial y_n} \\ 0\end{pmatrix}
\end{pmatrix}
\\&
= (x_1,\dots,x_n,1)^t\otimes (\tfrac{\partial}{\partial y_1},\dots,\tfrac{\partial}{\partial y_n},0) +  
 (\tfrac{\partial}{\partial y_1},\dots,\tfrac{\partial}{\partial y_n},0)^t\otimes (x_1,\dots,x_n,1)
\end{align*}
This has rank 2 and so the determinant of each $3\x 3$-minor vanishes.
\hfill $\square$


\subsection*{Schwartz space on $\cAL(n,\R)$.}
For an element $X$ of the Lie algebra of $\ASp(2n,\R)$ denote by
$L_X$ the corresponding vector field on the space of functions on $\cAL(n,\R)$.
We say that a $C^\infty$-function $f$ on $\cAL(n,\R)$ is contained in
the space $\cS(\cAL(n,\R))$ if for any collection $X_1$, \dots $X_k$
and each $M>0$
we have an estimate
$$
|L_{X_1}\dots L_{X_k} f(P)|=O(\text{distance between $0$ and $P$})^{-M}
.
$$ 

\begin{lemma}
A $C^\infty$ function $f$ is contained in $\cS(\cAL(n,\R))$ if and only if
for any chart of our atlas and for any compact subset $\Omega\subset\Symm$ 
any partial derivative of $f(T,\tau)$ of any order
is $O(\tau)^{-M}$ for all $M$.
\end{lemma}

{\sc Proof.}
Denote the functional space satisfying the condition of the lemma by $\cS(\cAL(n,\R))^{bis}$.

In coordinates, the vector fields $L_X$ are first order differential operators  whose coefficients are 
polynomials of degree $\le 2$ (formulas can be easily obtained from  (\ref{eq:perenos-1})).
This implies $\cS(\cAL(n,\R))\supset \cS(\cAL(n,\R))^{bis}$.
  Consider the subgroup in $\ASp(2n,\R)$ consisting of matrices
$$
\begin{pmatrix}
 1&b&\beta\\
0&1&0\\
0&0&1
\end{pmatrix}
$$
The Lie algebra  action of this subgroup is spanned
by the operators $\frac\partial{\partial t_{ij}}$ and 
$\frac\partial{\partial \tau_j}$. This implies $\cS(\cAL(n,\R))\subset \cS(\cAL(n,\R))^{bis}$.
\hfill $\square$

\begin{lemma}
The image of $\cS(\R^{2n})$ under the Radon transform is contained in  
$\cS(\cAL(n,\R))$.
\hfill $\square$
\end{lemma}



\subsection*{The $\cR_{flat}$-image of the Schwartz space.}

\begin{theorem}
\label{th:image}
A function $f$ on $\Symm(n,\R)\times\R^n$ is contained in the image
of $\cR_{flat}$ if and only if $f$ is contained in $\cS(\cAL(n,\R))$ and 
satisfies the determinantal system {\rm(\ref{eq:system})--(\ref{eq:system-2})} of PDE.
\end{theorem}

The proof occupies the next section.

\subsection*{Equations for $\cR_{comp}f$}
Now let $g=\cR_{comp}f$ be real analytic.
Consider its holomorphic extension $G$ to a small neighborhood $\cO$
of the manifold consisting of all $(S,\sigma)$ with
$$
SS^*=\bbI,\quad S=S^t, \quad \ov\sigma=- S^{-1}\sigma
$$
in $\Symm(n,\C)\times \C^n$. The system $\cD_{flat}$ transforms to a 
system $\cD_{comp}$ in a neighborhood of this manifold.  The following statement describes this system $\cD_{comp}$.

\begin{theorem}
\label{th:comp-determinantal}
Let $G$ be such a holomorphic extension. Then
the  {\rm(}two-valued{\rm)} analytic  function
\begin{equation}
\det(S)^{-1/2} G(S,\sigma)
\label{eq:system-twisted}
\end{equation}
satisfies the determinantal system $\cD_{flat}$, see {\rm (\ref{eq:system})--(\ref{eq:system-2})}.
\end{theorem}

\subsection*{Invariance of the determinantal system}

\begin{theorem}
\label{th:invariance}
{\rm a)}
Let $\begin{pmatrix}a&b&r\\
 c&d&s\\
0&0&1
\end{pmatrix}\in\ASp(2n,\C)$.
Let $\cU\subset \Symm(n,\C)\times\C^n$ be a set such that
$\det(cT+d)\ne 0$ for $(T,\tau)\in\cU$.
If a holomorphic function $G(T,\tau)$ defined  on $\cU$ satisfies the determinantal system 
{\rm(\ref{eq:system})--(\ref{eq:system-2})},
then the function
\begin{equation}
G\bigl((aT+b)(cT+d)^{-1},\,\,  a\tau+r\, -(aT+b)(cT+d)^{-1} (c\tau+s)  \bigr)\cdot \det(cT+d)^{-1}
\label{eq:transformation}
\end{equation}
also satisfies the determinantal system.

\sm

{\rm b)}
The same statement holds for smooth solutions of the determinantal system
defined on open subsets on $\Symm(n,\R)\times\R^n$. Such solutions also are invariant with respect to transformations {\rm(\ref{eq:perenos-1})}.

\sm

{\rm c)} Let $g$ be a real analytic function determined on an open subset
of $\USymm(n)\times\R^n$. Let its holomorphic continuation satisfies
the determinantal system $\cD_{comp}$.
 Then holomorphic continustion of {\rm(\ref{eq:perenos-2})} satisfies the same system.
\end{theorem}

Theorems \ref{th:comp-determinantal} and
\ref{th:invariance} are proved in Section \ref{sec:ImSchwartz}.


\section{\label{sec:ImSchwartz}Image of the Schwartz space. Proof of Theorem \ref{th:image}}

\subsection*{Fiber-wise Fourier transform on $\cAL(n,\R))$.%
\label{ss:fiber-wise}}
Let us fix an Euclidean inner product on $\R^{2n}$. For a linear subspace $Y\subset \R^{2n}$ denote
by $Y^\bot$ its orthocomplement with respect to  the Euclidean inner product.
For a vector $\xi$   we denote by $\xi^\bot$ the orthocomplement to the line $\R\xi$.

For a function $F\in\cS(\cAL(n,\R))$ we define its fiber-wise Fourier transform
$\cH_{comp} F$ on $\cS(\cAL(n,\R))$ in the following way. Let $L$ be a Lagrangian {\it linear} subspace.
Let
$$
\xi\bot L.
$$
Then
$$
\cH_{comp} F(\xi+L)=\int_{L^\bot} F(\eta+L) e^{i\la \xi,\eta\ra}\,d\eta.
$$

\begin{lemma}
If $F\in\cS(\cAL(n,\R))$ then $\cH_{comp} F\in\cS(\cAL(n,\R))$.
\hfill \qed
\end{lemma}

\begin{proposition}
\label{pr:fiber-wise}
Let $F= \cR_{comp} f$, where $f\in\cS(\R^{2n})$. For fixed $\xi$ and $L$
orthogonal to $\xi$, the function $\cH_{comp} F(\xi+L)$ is independent of $L$.
The function 
$$\phi(\xi):=\cH_{comp} F(\xi+L)$$
coincides with the Fourier transform $\wh f$  of $f$
\end{proposition}

{\sc Proof.} Let $\xi\in\R^{2n}$, and let $L$ be a linear Lagrangian 
subspace contained in $\xi^\bot$. We have:
\begin{multline*}
\wh f(\xi)=\int_{\R^{2n}}f(x)e^{i\la\xi, x\ra}\,dx
=\int_\R e^{is \la\xi,\xi\ra} \int_{z\in\xi^\bot} f(s\xi+z) dz\,ds
=\\=
\int_\R e^{is \la\xi,\xi\ra} \int_{v\in {(\R\xi+L)}^\bot}\int_{u\in L} f(s\xi+u+v) du\,dv\,ds
=\\=
\int_\R e^{is \la\xi,\xi\ra} \int_{v\in {(\R\xi+L)}^\bot} \cR f(s\xi+v+L)\,dv\,ds
=\int_{\eta\in L^\bot} e^{i\la \xi,\eta \ra} \cR f(\eta+L)\,d\eta.
\end{multline*}
So we get the desired statement.
\hfill $\square$

\subsection*{Reformulations.}
Now we describe the same transform in invariant language.

Denote by $\cAL^\circ(n,\R)$ the set whose points are pairs
$(L,\ell)$, where $L$ is a
Lagrangian linear subspace and $\ell$ is a linear functional on $\R^{2n}$ annihilating
$L$. 

{\sc Remark.}
Notice that a Euclidean inner product on $\R^{2n}$ identifies $\cAL(n,\R)$ and $\cAL^\circ(n,\R)$
(but this identification depends on a choice of the inner product). On the other hand 
the skew-symmetric bilinear form $\Lambda$ provides a canonical isomorphism between
$L$ and $\R^{2n}/L$.

\sm

For a function $F$ on $\wt \cAL(n,\R)$ satisfying the homogeneity
conditions (\ref{eq:homogen}) we will construct a function
$\cH F$ on $\cAL^\circ$.

Fix a linear Lagrangian subspace 
$L\subset\R^{2n}$ and a volume form $\omega$ on $L$. There is a well-defined form
$\omega^\circ$ on $\R^{2n}/L$ such that 
$dx\,dy=\omega\wedge \omega^\circ$. 
Then
$$
\cH_{comp} F(L,\ell)=\int_{\eta\in \R^{2n}/L} e^{i\ell(\eta)} F(\eta+L,\omega)d\omega^\circ
.
$$
The result does not depend on the choice of $\omega$.

\smallskip

We also write the same transform in the flat coordinates
$(T,\tau)\in\Symm(n,\R)\times \R^n$:
$$
\cH_{flat} F (T,v)=\int_{\R^{n}} F(T,\tau) e^{i v^t \tau}\,d\tau
.
$$
It is simply a Fourier transform with respect to a part of coordinates.

\subsection*{Inverse statement.}

\begin{proposition}\label{prop:inverse}
Let $F\in\cS(\wt\cAL(n,\R))$ satisfy the determinantal system.
 Then for any $\xi\ne 0$ and $L$ orthogonal to $\xi\ne 0$ the function
$$\frf(\xi):=\cH_{flat} F(\xi+L)$$
 is independent of $L$.
\end{proposition}

\begin{lemma}\label{lem:determinantal}
 Let $F$ satisfy the determinantal system $\cR_{flat}$. Then 
the function $g(T,v)=\cH_{flat} F(T,v)$
satisfies the condition
\begin{equation}
g(T+S, v)= g(T,v)\qquad    \text{if $v\ne 0$, $Sv=0$}
\label{eq:Sv}
\end{equation}
Equivalently, $g$ depends only on $(Tv,v)$.
\end{lemma}

Thus we get a function $\frf(w,v)$ on $\R^{2n}$, $w\ne 0$,  determined by
$$
\frf(Tv,v)=g(T,v)
.
$$

{\sc Proof.} We apply the transform $\cH_{flat}$ to the equation
$$
\det\begin{pmatrix}
2\frac{\partial}{\partial t_{kk}}& \frac{\partial}{\partial t_{kl}}& \frac{\partial}{\partial \tau_k}\\
\frac{\partial}{\partial t_{kl}}& 2\frac{\partial}{\partial t_{ll}}& \frac{\partial}{\partial \tau_l}\\
\frac{\partial}{\partial \tau_k}& \frac{\partial}{\partial \tau_l}&0
\end{pmatrix}
F(T,\tau)=0
$$
and obtain
\begin{equation}
\left(
v_k^2 \frac{\partial}{\partial t_{ll}} - v_k v_l \frac\partial{\partial t_{kl}}
+v_l^2 \frac{\partial}{\partial t_{kk}}
\right) g(T,v)=0
\label{eq:field}
.
\end{equation}
We interpret \eqref{eq:field} as follows: The function $g(T,v)$
is annihilated by a vector field.
The corresponding flow on $\Symm(n,\R)\times \R^n$ of this vector field is given by
$(T,v)\mapsto (T+\epsilon S^{kl},v)$, where 
$S^{kl}$ is the matrix, whose nonzero entries are
$$
s_{kk}=v_l^2,\qquad s_{kl}=s_{lk}=-v_kv_l,\qquad s_{ll}=v_k^2.
$$
It satisfies $S^{kl}v=0$.
 
Let all entries of $v$ be non-zero.
Any symmetric matrix $R$ satisfying $Rv=0$ is a linear combination
of matrices $S^{kl}$. Indeed,
we have
$$
\left(R-\sum_{k>l} \frac{r_{kl}}{v_k v_l} S^{kl}\right)v=0
$$
 The matrix in brackets is diagonal. Therefore it is zero.
 
 Thus, for  $v$  with non-zero entries we have (\ref{eq:Sv}). 
 Next, our system (\ref{eq:field}) is invariant with respect to the transformations
 $T\mapsto ATA^t$, $v\mapsto Av$. Therefore we can send each nonzero vector $v$ to a vector
 with nonzero entries.
 \hfill
 $\square$

Now consider a function $F\in \cS(\cAL(n,\R))$ such that 
$F(\xi+L)$ is independent of $L$ if $L\bot \xi$, $\xi\ne 0$.
Then we write a function $\frf(\xi):=F(\xi+L)$.

\sm
  
{\sc Proof of proposition \ref{prop:inverse}.}
Indeed, let $\xi\ne 0$, $L\in \cL(n,\R)$.
We chose a pair of transversal linear Lagrangian subspaces $V_-$, $V_+\subset \R^{2n}$
in such way that $L\cap V_+=0$. Then we introduce coordinates $(T,\tau)$
and apply the lemma.
\hfill $\square$

\subsection*{Smoothness at zero.}  
  
\begin{proposition}
The function 
$\frf$ admits a smooth continuation at the point $\xi=0$.

\end{proposition}

{\sc Proof.}
 The condition (\ref{eq:Sv}) is equivalent to 
(\ref{eq:field}).

We introduce a $\Z^n$-gradation in the space of polynomials in the variables $t_{kl}$ and 
$v_m$. Denote by $\delta_j$ the natural basis of $\Z^n$.
We set     
$$
\deg  t_{kl}=\delta_k+\delta_l,\qquad \deg v_l=\delta_l.
$$
The differential operators
(\ref{eq:field}) are homogeneous with respect to this gradation.
Consider the Taylor expansion of $f(T,v)$ at zero,
$$
f(T,v)=\sum_{m=0}^N f^{(m)}(T,v)+ R_N(T,v)
.
$$
By homogeneity of the differential operators
 each polynomial $f^{(m)}$ satisfies the equations
 (\ref{eq:field}). Therefore the remainder $R_m(T,v)$ satisfies the same equations. Thus all 
 $f^{(m)}$ and $R_N$ depend only on $(Tv,v)$, by Lemma \ref{lem:determinantal}
 
For each $m$, we consider the extension of the polynomial $f^{(m)}$ to  the complex linear space  $\Symm(n,\C)\times \C^n$.
We pass over to the function 
$$
\frf^{(m)}(w,v)=f^{(m)}(Tv,v)
,$$
which is holomorphic outside the set $v=0$. The codimension of this set is 
$n>1$. Therefore by Hartogs' theorem
the function $\frf^{(m)}$ is holomorphic on the whole of $\C^n$.

We consider now the remainder
$$
R_N(T,v)=O(\sum |t_{ij}|+ \sum |v_j|)^{M+1},\qquad T,v\to 0
$$
 Consider
the function 
$$
\frR(w,v)=R^{(m)}(Tv,v)\qquad\text{if $w=Tv$.}
$$
Notice that if $|T|\le \epsilon$, $|v|\le \epsilon$, then
$|Tv|\le \epsilon^2$, on the hand any vector $w$ such that $|v|<\epsilon^2$
can be represented in such form.
Therefore
$$\frR(w,v)=O(\sum |w_j|+\sum|v_j|)^{M/2},\qquad v\ne 0.$$
 
 Hence $\frf$ is $C^\infty$ at zero.
 
\subsection*{End of proof of Theorem \ref{th:image}.}
Let $F\in\cS(\cAL(n,\R))$, $\frf\in\cS(\R^{2n})$
be the same as above.

\begin{lemma}
Let 
$$
g(x)=\int_{\R^{2n}}
\frf(\xi)e^{-i\xi^t x}\,dx
$$
Then 
$$
F(\eta+L)=\cR_{comp} g(\eta+L), \qquad \eta\bot L
$$
\end{lemma}

{\sc Proof.}  
We can evaluate the Fourier transform $f(\xi)$ of a function $g\in\cS(\R^{2n})$
in the following way. For $\xi\in\R^{2n}$ consider a Lagrangian subspace 
$L$ orthogonal to $\xi$. 
Then
$$
f(\xi)=\int_{L^\bot} \cR_{comp}(x+L)e^{i\xi^t x}\,dx
.$$

\section{Invariance of the determinantal system. Extraneous solutions}

\subsection*{Proof of Theorem \ref{th:invariance}. Holomorphic case}
Let the matrix $T=\{t_{kl}\}$ range in $\Symm(n,\C)$. 
Consider the following matrix
\begin{equation}
\begin{pmatrix}
2\frac \partial{\partial t_{11}}&\frac \partial{\partial t_{12}}&\dots\\
\frac \partial{\partial t_{21}}&2\frac \partial{\partial t_{22}}&\dots\\
\vdots&\vdots&\ddots
\end{pmatrix}
\label{eq:partial-matrix}
.
\end{equation}
We say that \emph{a function $F$ on $\Symm(n,\C)$ satisfies the determinantal system 
$\Delta_k$ if $F$ is annihilated by all $(k+1)\times (k+1)$-minors of this matrix.} 

\begin{theorem}
Let $\begin{pmatrix}
a&b\\c&d
\end{pmatrix}\in\Sp(2n,\C)
$. Let 
$\cO\subset \Symm(n,\C)$ be an open  subset, and $\det(cT+d)\ne 0$ on $\cO$.
Then for any holomorphic function $F$ on $\cO$ satisfying the determinantal system $\Delta_k$,
the function 
\begin{equation}
F\bigl((aT+b)(cT+d)^{-1}\bigr)\det(cT+d)^{-k/2}
\label{eq:book}
\end{equation}
also satisfies the determinantal system $\Delta_k$.
\end{theorem}

The statement was formulated in \cite{Neretin11}, Subsection 8.8.3, 
without formal proof, but modulo considerations of \cite{Neretin11}, Section 8.8,
a proof is straightforward. 

Set $k=2$.

Now consider the group $\Sp(2n+2,\C)$ acting (locally) on
$\Symm(n+1,\C)$ by the formula
$$
\rho\begin{pmatrix}\wt a&\wt b\\  \wt c&\wt d \end{pmatrix} f(S)=
F\bigl((\wt aS+\wt b)(\wt cS+\wt d)^{-1})\det(\wt cS+\wt d)^{-1}
$$
Restrict this action to the subgroup $\Sp(2n,\C)\subset \Sp(2n+2,\C)$.
Representing  elements of $\Symm(n+1,\C)$
as $(n+1)\times(n+1)$-matrices 
$
\begin{pmatrix}
S_{11}&S_{12}\\
S_{21}&S_{22}
\end{pmatrix}
$,
we get the following formula for the  action of $\Sp(2n,\C)$ 
  \begin{multline*}
 \rho \begin{pmatrix}
  a&0&b&0\\
  0&1&0&0\\
  c&0&d&0\\
  0&0&0&1
  \end{pmatrix}
  F
\begin{pmatrix}
S_{11}&S_{12}\\
S_{21}&S_{22}
\end{pmatrix}
  = \\
F\begin{pmatrix}
(aS_{11}+b)(cS_{11}+d)^{-1}
&
\bigl(a-(aS_{11}+b)(cS_{11}+d)^{-1}c\bigr)S_{12}
\\
S_{21}(cS_{11}+d)^{-1}& S_{22}-S_{21}(cS_{11}+d)^{-1}cS_{12})
\end{pmatrix}
\det (aS_{11}+b)^{-1}
\end{multline*}
Recall that $S_{22}$ is a number, $S_{12}$ is a column, and $S_{21}=S_{12}^t$.

We observe that:
\begin{itemize}
	\item 
The transformations $\rho$ send functions independent of $S_{22}$ to functions
independent of $S_{22}$.
\item
The action of $\Sp(2n,\C)$ on functions depending on $S_{11}$, $S_{12}$
coincides with the action (\ref{eq:transformation})
(we set $T=S_{11}$, $\tau=S_{12}$; since we work with
$\Sp(2n,\C)$, we have $r=0$, $s=0$).
\item
The restriction of the system $\Delta_2$ to functions independent of $S_{22}$ is precisely
the determinantal system (\ref{eq:system})--(\ref{eq:system-2}). 
\end{itemize}

Thus we proved the invariance of the determinantal system 
with respect to $\Sp(2n,\C)$. The invariance with respect
to translations is obvious.

\subsection*{Proof of Theorem \ref{th:invariance}. Smooth case}
Let $F$ be a function defined on an open subset in $\Symm(n,\R)\times\R^n$.
Consider the Taylor decomposition of $F$ at a point $\frs=(S,\sigma)$, 
$$
F(T,\tau)=\sum_{j=0}^m p^{j}_{S,\sigma} (T,\tau)+ R^m_{S,\sigma}(T,\tau)
,$$
where $p^{j}_{S,\sigma} (T,\tau)$ is a polynomial of degree $j$
homogeneous with respect to $t_{kl}-s_{kl}$, $\tau_k-\sigma_k$,
and $R^m_{S,\sigma}(T,\tau)$ has zero partial derivatives of order $\le m$ at
the point
$(S,\sigma)$

\begin{lemma}
A function $F$ satisfies the determinantal system
if and only if for all $(S,\sigma)$ for all $j$ polynomials $p^j_{S,\sigma} (T,\tau)$
satisfy the same conditions.
\hfill\qed
\end{lemma}


\sm

Let us prove the theorem. Obviously the determinantal system is invariant with respect to
translations $(T,\tau)\mapsto (T+H, \tau+h)$ and transformations
$$
(T,\tau)\mapsto (aTa^t, a\tau),\qquad a\in\GL(n,\R).
$$
In other words, the determinantal system is invariant with respect to the
parabolic subgroup $P\subset\Sp(2n,\R)$ consisting 
of matrices
\begin{equation}
\begin{pmatrix}
a&b&r\\
0&(a^t)^{-1}&s\\
0&0&1
\end{pmatrix}\in\ASp(2n,\R)
\label{eq:1}
\end{equation}
Let us add the subgroup $N_-\subset\Sp(2n,\R)$ consisting of matrices
\begin{equation}
\begin{pmatrix}
\bbI&0&0\\
c&\bbI&0\\
0&0&1
\end{pmatrix}\in\ASp(2n,\R)
;
\label{eq:2}
\end{equation}
it acts by transformations
\begin{equation}
F(T,\tau)\mapsto F\bigl(T(\bbI+cT)^{-1}, \big((\bbI+cT)^t\big)^{-1}\tau\bigr)\det(\bbI+cT)^{-1}
.
\label{eq:3}
\end{equation}

\begin{lemma}
If $F$ satisfies the determinantal system, then for any
minor $\Xi$ of {\rm(\ref{eq:system})} the right-hand side of {\rm(\ref{eq:3})} is flat
at 0.
\end{lemma}

{\sc Proof.} The Taylor polynomial of order $m$ at $(T,\tau)=(0,0)$ of the transformed function
depends only on the Taylor polynomial of order $m$ of $F$. On the other hand,
$\Xi$ applied to the Taylor polynomial is zero.
\hfill $\square$

\sm

Now take a point $(S_0,\sigma_0)$ and element $g$ of $\ASp(2n,\R)$.
First, we consider a translation $r$ sending $(S_0,\sigma_0)$ to $(0,0)$.
Next, decompose $gr^{-1}$ as a product $p\cdot h$, where $p\in P$ and
$h\in N$. Clearly, applying a minor $\Xi$ to the transformed function, we get a function 
which is flat at $g(S_0,\sigma_0)$.
\hfill $\square$

\subsection*{Modifications of formulas for transformations}
Thus we proved that the determinantal system $\Gamma$ is invariant with respect 
to transformations (\ref{eq:transformation}). But expressions 
(\ref{eq:perenos-1}), (\ref{eq:perenos-2}), (\ref{eq:perenos-3})
have a slightly different form, they 
contain absolute values of determinants.

1) In the case of  (\ref{eq:perenos-1}),
the expression 
$$\frac{|\det(cT+d)|}{\det(cT+d)}$$
is locally constant.

\sm

2)
Let us examine the transformations  (\ref{eq:perenos-2}).
We write
\begin{multline}
\bigl|\det (\ov \Psi S+\ov\Phi  )\bigr|^{-1}
=\\=
\det (\Phi\ov\Phi)^{-1/2} \cdot
\det(\bbI+ \Phi^{-1}\Psi S)^{-1/2} 
\cdot
\det(\bbI+ \ov\Phi^{-1}\ov\Psi S)^{-1/2}. 
\label{eq:PhiPhi}
\end{multline}
First, we note that
$$
\det (\Phi\ov\Phi)^{-1/2}=\det (\Phi \Phi^*)^{-1/2}
$$
and we can take the positive value of the square root ($\Phi$ is non-degenerate).
Next we evaluate $\det(\cdot)^{-1/2}$
as
$$
\det(\bbI+H)^{-1/2}:=\det (\bbI-\frac 12 H+\frac 18 H^2+\dots)
$$
and $\|H\|$ in our case is $\le 1$ (see, e.g., \cite{Neretin11}, Theorem 3.3.5).
Thus, for fixed 
$\begin{pmatrix}\Phi&\Psi&p\\
 \ov\Psi&\ov\Phi&\ov p
 \\
 0&0&1\end{pmatrix}$,
the expression (\ref{eq:PhiPhi}) as a function of $S$ is  well-defined and holomorphic  in a neighborhood
of $\U(n)\cap\Symm(n,\C)$ in  $\Symm(n,\C)$.

We also transform the same expression as 
\begin{multline*}
\bigl|\det (\ov \Psi S+\ov\Phi)\bigr|^{-1}
=\pm\det (\ov \Psi S+\ov\Phi)^{-1/2}\det  (\Psi S^{-1}+\Phi)^{-1/2}
=\\=
\pm
\det(\ov \Psi S+\ov\Phi)^{-1} \cdot
\left(\frac{\det\left((\Phi S+\Psi) (\ov \Psi S+\ov\Phi)^{-1}\right)}
{\det S}\right)^{-1/2} 
.
\end{multline*}
Therefore, we can write  (\ref{eq:perenos-2})
as
\begin{multline}
\rho_{comp}
\begin{pmatrix}\Phi&\Psi&p\\
 \ov\Psi&\ov\Phi&\ov p\\
0&0&1
\end{pmatrix}
\Bigl[f(S,\sigma)(\det S)^{1/2}\Bigr]=
\\
=
\pm
\biggl[
f\Bigl((\Phi S+\Psi)(\ov \Psi S +\ov\Phi)^{-1},  \Phi\sigma+p\, -(\Phi S+\Psi)(\ov\Psi S+\ov\Phi)^{-1} (\ov\Psi\sigma+\ov p)  \Bigr)
\cdot\\ \cdot
\det\left((\Phi S+\Psi) (\ov \Psi S+\ov\Phi)^{-1}\right)^{-1/2}\biggr]
\cdot
 \det (\ov \Psi S+\ov\Phi)^{-1}
\label{eq:perenos-5}
\end{multline}

Thus, {\it the twisted function
$$f(S,\sigma)(\det S)^{1/2}$$
transforms
according to formula} (\ref{eq:transformation}).

\sm

3) It remains to examine the transform (\ref{eq:perenos-3})
from the flat to the compact model. Again, we transform the last factor,
\begin{multline*}
|\det(\bbI+iS)|^{-1}
= \det(\bbI+iS)^{-1/2}\det(1-iS^{-1})^{-1/2}
=\pm\det (iS)^{1/2} \det(\bbI+iS)^{-1}.
\end{multline*}

Therefore formula  (\ref{eq:perenos-3})
takes the form
\begin{multline}
\cR_{comp}(S,\sigma)=\\
\cR_{flat}\bigl((S+i\cdot\bbI)(iS+\bbI)^{-1}, (\bbI+iS^{-1})^{-1}\sigma\bigr)
\cdot(\det S)^{1/2}\det(\bbI+iS)^{-1}\cdot 2^n.
\label{eq:perenos-4}
\end{multline}

\subsection*{Proof of Theorem \ref{th:comp-determinantal}}
It follows from formula (\ref{eq:perenos-4}).

\subsection*{Extraneous solutions}
All solutions of the determinantal system $\Gamma$ in the function space $\cS(\cAL(n,\R))$
are contained in the image of the Radon transform.
We wish to show that there are extraneous  solutions of the system defined on open subsets of
$\cS(\cAL(n,\R))$. 

Any function on $\cL\simeq\U(n)/\O(n)$ determines a function on $\cAL(n,\R)$, which is constant
on fibers. We will present extraneous solutions of the determinantal system $\cD_{comp}$ on the 
double covering 
$\U(n)/\O(n)^\sim$
of $\U(n)/\O(n)$.

The group $\U(n)$ acts on $\cL(n)$ by transformations
$u:\,T\mapsto uTu^t$.
Recall that all irreducible representations $T[m]$ of
$\U(n)$ are enumerated by integer signatures
$$
m:\,m_1\ge\dots\ge m_n
.
$$
 The space of $C^\infty$-functions on
$\cL(n,\R)$ is  a multiplicity-free direct sum
$$
\bigoplus_{p_1\ge\dots\ge p_n} T[2p_1,\dots,2p_n]
$$
(see, e.g. \cite{Helgason94}, \cite{Hua}).
The functions 
$$
\prod_{j=1}^n \det\begin{pmatrix} t_{11}&\dots& t_{1j}\\
\vdots&\ddots& \vdots\\
t_{1j}&\dots &z_{jj} \end{pmatrix}^{p_{j}-p_{j+1}}
$$
are highest weight vectors of these representations.

The space of smooth functions on $\U(n)/\O(n)^\sim$ 
is a direct sum
$$
\Bigl(\bigoplus_{p_1\ge\dots\ge p_n} T[2p_1,\dots,2p_n]\Bigr)
\oplus \Bigl(\bigoplus_{q_1\ge\dots\ge q_n} T[2q_1+1,\dots,2q_n+1]\Bigr)
.$$

We look for real analytic solutions $f$ of the  determinantal system 
$\cR_{comp}$
 on  $\U(n)/\O(n)^\sim$. 
 This is equivalent to the determinantal system $\Delta_3$ applied to functions
 $$\det(T)^{-1/2}f(T),$$
 where $f$ is analytic in a neighborhood of the surface
 $TT^*=1$.

\begin{proposition}
\label{pr:extraneous}
All representations of $\U(n)$ with the following signatures 
are contained in the space of solutions of $\Delta_3$ on $\U(n)/\O(n)^\sim$
\begin{align}
&(2\alpha+1,2\beta+1,1,\dots,1),\qquad & \alpha\ge \beta\ge0
\label{eq:signature-1}
\\
&(2\alpha,0,\dots, 0, -2\gamma),\qquad & \alpha\ge 0, \gamma\ge 0
\label{eq:signature-2}
\\
&(-1,\dots,-1, -2\gamma-1,-2\delta-1),
\qquad & \delta\ge\gamma \ge 0
\label{eq:signature-3}
\end{align}
\end{proposition}

Observe that solutions of the type (\ref{eq:signature-2})
are global solutions on $\cL(n,\R)$,
and solutions 
 (\ref{eq:signature-1}), (\ref{eq:signature-3}) are determined on the double covering.

{\sc Proof.}
Functions 
$$t_{11}^{\alpha-\beta}\det\begin{pmatrix}t_{11}&t_{12}\\ t_{12}&t_{22}\end{pmatrix}^\beta$$
obviously are annihilated by all $3\times 3$ minors of the 
matrix (\ref{eq:system}). By $\U(n)$-invariance we get
that all operators of $\Delta_3$ anihilate all subspaces 
$T[2\alpha,2\beta,0,\dots,0]$.
Multiplication by $\det(T)^{1/2}$ 
implies a shift of signatures 
by $(1,\dots,1)$.
Thus we get signatures of the form (\ref{eq:signature-1}).

We apply the transformation
$$
\pi_2\begin{pmatrix}
0&\bbI\\
-\bbI&0
\end{pmatrix} g(z)=g(T^{-1})\det(T)^{-1}
$$
(it is of the  form (\ref{eq:book}) with $k=2$)
and obtain that the functions 
$$
\det
\begin{pmatrix} t_{11}&\dots& t_{1(n-2)}\\
\vdots&\ddots& \vdots\\
t_{1(n-2)}&\dots& t_{(n-2)(n-2)} \end{pmatrix}^\beta
\det
\begin{pmatrix} t_{11}&\dots& t_{1(n-1)}\\
\vdots&\ddots& \vdots\\
t_{1(n-1)}&\dots& t_{(n-1)(n-1)} \end{pmatrix}^{\alpha-\beta}
\det
T^{-\alpha-1}.
$$
are also solutions of $\Delta_3$.
This gives signatures of the type (\ref{eq:signature-3}).

Next, the functions
\begin{equation}
t_{11}^\alpha
\label{eq:mnozh-1}
\end{equation}
are annihilated by $2\times 2$ minors of the 
matrix (\ref{eq:partial-matrix}). We apply the transformations
$$
\pi_1 f(T)\begin{pmatrix}
0&\bbI\\
-\bbI&0
\end{pmatrix} g(T)=g(T^{-1})\det(T)^{-1/2}
$$
and observe that the functions
\begin{equation}
\det
\begin{pmatrix} t_{11}&\dots& t_{1(n-2)}\\
\vdots&\ddots &\vdots\\
t_{1(n-2)}&\dots& t_{(n-2)(n-2)} \end{pmatrix}^\alpha \det T^{-\alpha-1/2}
,
\label{eq:mnozh-2}
\end{equation}
are also annihilated by $2\times 2$ minors of the matrix
(\ref{eq:partial-matrix}).
Therefore products of functions 
(\ref{eq:mnozh-2})  (\ref{eq:mnozh-1}) are anihilated by 
$3\times 3$-minors of (\ref{eq:partial-matrix}).
Thus we get get signatures of the type (\ref{eq:signature-2}).
\hfill $\square$


\section{Inversion formulas}

We wish to reconstruct a function $f(x,y)$ of $2n$ variables from a function
$F=\cR f$ of $n(n+1)/2+n$
variables. There are many ways to do this, we propose 3 variants
(which are straightforward imitations of inversion formulas for the usual Radon transform).

\subsection*{Inversion 1}
As we have seen in Proposition  \ref{pr:fiber-wise}, the Fourier transform of
$f$ can be reconstructed from $\cR f$.

\subsection*{Inversion 2}
Introduce  parameters $\phi_1$, \dots, $\phi_n\in[0,2\pi)$, and
$p_1$, \dots, $p_n\ge 0$.
For each $\phi$, $p$
consider the    affine Lagrangian  subspace $L=L[\phi,p]$
 defined by the equations
$$
\begin{cases}
x_j \cos\phi_j+y_j\sin\phi_j=p_j, \quad j=1,\dots,n .
\end{cases}
$$
The set $\bbA(n)$ of all such subspaces is a submanifold $\cAL(n,\R)$.

Consider the decomposition
$$\R^{2n}=\R^2\times\dots\times\R^2,$$
where each factor has coordinates $(x_j,y_j)$.
In $j$-th factor consider a line $\ell[\phi_j,p_j]$  defined
by the equation    $x_j \cos\phi_j+y_j\sin\phi_j=p_j $.
The Lagrangian subspace $L[\phi,p]$ is the direct product of the lines
$\ell[\phi_j,p_j]$, moreover
\begin{align}
d\lambda_{L[\phi,p]}=d\lambda_{\ell[\phi_1,p_1]}\times\dots \times d\lambda_{\ell[\phi_n,p_n]};
\\
d\mu_{L[\phi,p]}=d\mu_{\ell[\phi_1,p_1]}\times\dots \times d\mu_{\ell[\phi_n,p_n]},
\end{align}
and
$$
\bbA(n)=\cAL(1)\times\dots\times \cAL(1)
$$

Restrict $\cR_{comp}f$ to $\bbA(n)$. We get the integral transformation
$$
\cR f(\phi,p)=
\int_{L[\phi,p]} f(x,y) \,d\mu_{L[\phi,p]}
$$
the left hand side is a function on the product of the $n$-dimensional torus $\bbT^n$ and $\R^n$.
This transform is a tensor product of two-dimensional Radon transforms.
Therefore we can automatically write the inverse operator
as a tensor product of inverse operators for $\R^2$ (see, e.g. \cite{GGG}, \cite{Helgason}.
Denote by $F$ the average of $\cR f$  over the torus
$$
F(p_1.\dots,p_n)=\frac 1{(2\pi)^n}\int_0^{2\pi} \dots \int_0^{2\pi}
\cR f(\phi_1,\dots,\phi_n, p_1,\dots,p_n)]\,d\phi_1\dots d \phi_n
$$
Then
$$
f(0,\dots,0)=\int_0^\infty \dots \int_0^\infty 
\frac 1 {p_1}\dots \frac 1 {p_n}
\frac \partial{\partial p_1}\dots\frac \partial{\partial p_1} F(p_1,\dots,p_n)
\,dp_1\dots dp_n
$$

The problem is translation invariant, therefore we can get the value of $f$ at any point. 

\subsection*{Inversion 3}
The usual inversion formula also can be applied to the
Lagrangian Radon transform.
Denote by $G(v)$ the average of $F=\cR_{comp}f$  over
all affine Lagrangian subspaces containing a point $v$.
 Denote by $\sigma_n$ the area  of $(k-1)$-dimensional unit sphere,
$$
\sigma_k=\frac{2\pi^{k/2}}{\Gamma(n/2)}
$$

\begin{lemma}
 \begin{equation}
G(v)=\frac{\sigma_n}{\sigma_{2n}}
\int_{\R^{2n}} f(w) |v-w|^{-n} \,dv
\label{eq:riesz}
\end{equation}
\end{lemma}

{\sc Proof.} Denote by $d\nu_L$ the measure
$d\mu_L$ regarded as a measure on the whole $\R^{2n}$.
We consider the average of all $d\nu_L$ for $L$ passing through
$v$ and integrate $f$ with respect to this measure
\hfill $\square$

\sm

The integral equation (\ref{eq:riesz}) is well known, it can be solved  in terms of the Riesz potentials,
see, e.g. \cite{Helgason}. Consider the operator
$$
I_\alpha f(v)=\frac{2^{-\alpha}\pi^{}\Gamma(n-\alpha/2)}{\Gamma(\alpha/2)}\int_{\R^{2n}} |v-w|^{-2n+\alpha}
f(w)\,dw
.$$
For a fixed $f\in\cS(\R^{2n})$, we get a function $I_\alpha f$ admitting a meromorphic
continuation in $\alpha$ to $\C$ with (possible) poles at
$\alpha\in 2n+ 2\Z_+$.
Then
$$
I_{-n} I_nf=f, \qquad f\in\cS(\R^{2n}),
$$
and this is the formula reconstructing $f$ from $G$ in (\ref{eq:riesz}).
Also, positive integer powers of the Laplace operator $\Delta$ can be expressed as
$$
\Delta^k=I_{-2k} 
$$
In particular, for even $n$, 
$$
f=\mathrm{const}\cdot\Delta^{n/2} G
$$

\section{Radon transform of distributions}

\subsection*{Counterexample.}
Consider the function $\gamma(x)$ on $\R^{2n}$ given by
$$
\gamma(x)=\frac 1{(1+|x|)^{n} \ln(2+|x|)}
$$
where $|x|$ denotes the distance between $0$ and $x$.
Let $h(r)$ be an arbitrary positive function on $\R^+$ such that
$$
h(r)=\begin{cases}
1,&\qquad r<1;
\\
0,& \qquad r >2.
\end{cases}
$$
The sequence
$$
\gamma_N(x):=h(|x|/N)\,\gamma(x)
$$
converges to $\gamma(x)$ in the space $\cS'(\R^{2n})$.
On the other hand the sequence
$$
\cR_{comp} \gamma_N (L)
$$
increases monotonically and  tends  to $\infty$
uniformly on compact sets. Therefore 
$\cR_{comp} \gamma_N (L)$ has no limit in $\cS'(\cAL(n,\R)$
(and also has no limit in the space $\cE'(\cA(L)$ of all distributions).

Notice that $\gamma\in L^2(\R^{2n})$. Therefore we can not extend the Radon transform
$\cS(\R^{2n})\to\cS(\cAL(n,\R))$
to a continuous operator $L^2(\R^{2n})\to \cE'(\cAL(n,\R))$.

\subsection*{Radon transform of distributions}
Denote by $\cB(\R^{2n})$ the set of all $C^\infty$ functions $f$ on
$\R^{2n}$ such that their partial derivatives
satisfy the estimates
$$
H_\alpha(f):=
\int_{\R^{2n}}
\sup (1+|x|)^{n} \left| \frac{\partial^\alpha}{\partial x^\alpha} f(x)\right| \,dx<\infty
$$
for any multiindex $\alpha=(\alpha_1,\dots,\alpha_{2n})$.
We regard $H_\alpha(f)$ as seminorms in $\cB(\R^{2n})$.
By $\cB'(\R^{2n})\subset \cS'(\R^{2n})$ we denote
the space of all distributions,  
which are continuous functionals on the space $\cB(\R^{2n})$.

By $\cAL_x$ denote the set of all $L\in\cAL(n,\R)$ containing $x$.
Consider the stabilizer $G_x\subset \U(n)\ltimes \R^{2n}$
of the point $x$. Denote by $d\sigma_{x}(L)$ the unique
probability  measure
on $\cAL_x$ invariant with respect to $G_x$.

\begin{theorem}
\label{th:distributions}
The Radon transform $\cR_{comp}:\cS(\R^{2n})\to\cS(\cAL(n,\R))$ admits 
a unique continuous extension to an operator
$\cB'(\R^{2n})\to\cS'(\cAL(n,\R))$.
For $\phi\in\cS(\cAL(n,\R))$ we 
define the function  $\fra \phi(x)$ on $\R^{2n}$ by 
$$
\fra\phi(x)=\int_{\cAL_x} \phi(L)\,d\sigma_x(L)
.$$
 For $f\in \cB'(\R^{2n})$ and $\phi\in\cS(\cAL(n,\R))$
 we set
$$
(\cR_{comp} f, \phi)=
(f, \fra\phi)
.
$$
\end{theorem}

The proof of Theorem \ref{th:distributions} will fill the rest of this section. 

\begin{lemma}
There is a 
$\U(n)\ltimes \R^{2n}$-invariant measure $d\sigma$ on 
the set $\cM$ of pairs $x\in \R^{2n}$, $L\in\cAL(n,\R)$ such that $x\in L$, which is unique up to a 
scalar factor. 
It is can be represented as 
$$
d\sigma(x,L)=d\sigma_x(L)\,dx= d\mu_L(x)\, dL
,$$
where $dL$ is the $\U(n)\ltimes \R^{2n}$-invariant measure on $\cAL(n,\R)$.
\end{lemma}

{\sc Proof.} The space $\cM$ is $\U(n)\ltimes \R^{2n}$-homogeneous, a stabilizer
of a point is isomorphic to $\O(n)$. Since the Haar measure on $\U(n)\ltimes \R^{2n}$
is two-side invariant and $\O(n)$ is compact, the quotient space $(\U(n)\ltimes \R^{2n})/\O(n)$
admits a unique invariant measure.
It is easy to verify that the measures $d\sigma_x(L)\,dx$ and  $d\mu_L(x)\, dL$
also are invariant.
\hfill $\square$

\begin{corollary}
For $f\in\cS(\R^{2n})$, $\phi\in \cAL(n,\R)$, we have
$$
\int_{\cAL(n,\R)} \cR_{comp} f(L)\phi(L)\,dL
=
\int_{\R^{2n}} f(x)\int_{\cAL_x} \phi(L)\, d\sigma_x(L)\,dx
.
$$
\end{corollary}

\begin{proposition}
\label{pr:decreasing}
If $\phi\in \cS(\cAL(2n))$, then
$$
\fra\phi(x)=\int_{\cAL_x} \phi(L)\, d\sigma_x(L)=O(|x|^{-n}),\qquad x\to\infty
.
$$
\end{proposition}

\subsection*{Estimate of $\fra\phi$.}

Denote by $\rB[x,r]\subset \R^{2n}$  the ball with center $x$ of radius $r$.

Denote by $\Gamma_r(x)$ the set of all subspaces $L\in \cAL_x$,
such that $\rB[0,r]\cap L\ne\varnothing$. Then
$\sigma_x(\Gamma_r(x))$ depends only on $\rho=r/|x|$; we
denote
$$
h(\rho):=\sigma_x(\Gamma_r(x)).
$$

\begin{lemma}
\label{l:asymptotics}
The function $h(\rho)$ monotonically increases, $h(0)=0$, $h(\rho)=1$ for
$\rho\ge 1$,
and
$$
h(\epsilon)=C\epsilon^n+ o(\epsilon^n),\qquad \epsilon \to 0
$$
\end{lemma}

{\sc Proof.}  Only the last statement requires a proof. Consider the product
$\rS^{2n-1}\times \cL(n,\R)$, where $\rS^{2n-1}$ is the unit sphere. Equip this space
with a natural probabilistic measure. Fix $x_0\in \rS^{2n-1}$ and $L_0\in\cL(n,\R)$.
Let:
\begin{itemize}
  \item 
$\Sigma\subset \rS^{2n-1}\times \cL(n,\R)$ be the set of all pairs
$x\in \rS^{2n-1}$, $L\in\cL(n)$ such that 
$L\cap \rB[x,r]\ne\varnothing$;
\item
$\Delta\subset \cL(n,\R)$ be the set of all $L\in\cL(n,\R)$ such that 
 $L\cap \rB[x_0,r]\ne\varnothing$;
\item
$\Xi\subset \rS^{2n-1}$ be the set of all $x\in \rS^{2n-1}$
such that $L_0\cap \rB[x,r]\ne\emptyset$.
\end{itemize}

Then 
$$\bigl\{\text
{measure of $\Delta$}\bigr\}
= 
\bigl\{\text
{measure of $\Sigma$}\bigr\}
=
\bigl\{\text
{measure of $\Xi$}\bigr\}
$$  
Let us equip the unit sphere $\rS^{2n-1}$ with the intrinsic metric.
The set $\Xi\subset \rS^{2n-1}$ is a  neighborhood of an $(n-1)$-dimensional equator of
radius $\arcsin r$. Its volume can be evaluated explicitly, but
asymptotics of volume  as $r\to 0$ is obvious.
\hfill $\square$

\sm

Denote by $\Theta_r$ the function on $\cAL(n,\R)$
given by
$$
\Theta_r(L)=\begin{cases}
1,\qquad\text{if $L\cap \rB[0,r]\ne \varnothing$};
\\
0, \qquad \text{if $L\cap \rB[0,r]= \varnothing$}.
\end{cases}
$$
Then 
$$
\fra \Theta_r(x)=h\bigl(r/|x|\bigr)
.
$$


\begin{lemma}
For $\phi\in \cS(\cAL(n,\R))$,
$$
\fra\phi(x)=O(|x|^{-n}),\qquad x\to\infty
.
$$
\end{lemma}

{\sc Proof.}
Denote by $d(0,L)$ the distance between $0$
and a subspace $L\in\cAL(n,\R)$. It can be easily shown
that $d(0,L)^2$ is a smooth function on $\cAL(n,\R)$.
It is sufficient to prove our statement
for functions $\phi(L)$ having the form
$$
\phi(L)=q(d(0,L))
,$$
where $q(t)\in\cS(\R)$ is an
even and monotonically decreasing 
on $[0,\infty)$.

We represent $\phi(L)$ as
$$
\phi(L)=\int_0^\infty \Theta_r(L)(- q'(r))\,dr.
$$
and get
$$
\fra\phi(x)=
\int_0^{\infty} \fra\Theta_r(x)(- q'(r))\,dr
=\int_0^{|x|} \fra\Theta_r(x)(- q'(r))\,dr + \int_{|x|}^\infty\fra\Theta_r(x) (-q'(r))\,dr
.$$
We have $\Theta_r(x)=1$ if $|x|<r$. Therefore the second summand
is 
$$
\int_{|x|}^\infty(- q'(r))\,dr= q(|x|)
.$$
Next, we take $\delta$ and C such that $h(\rho)\le C \rho^{n}$ for
$0<\rho<\delta$ and split $\bigl[0,|x|\bigr]$ into two segments
$\bigl[0,\delta|x|\bigr]$ and $\bigl[\delta|x|,|x|\bigr]$
\begin{multline*}
\int_0^{\delta|x|} \fra\Theta_r(x)(- q'(r))\,dr
=\int_0^{\delta|x|} h(\delta/|x|)(- q'(r))\,dr
\le\\ \le C |x|^{-n}
\int_0^{\delta|x|} r^n (-q'(r))\,dr
\le  C |x|^{-n}
\int_0^{\infty} r^n (-q'(r))\,dr\sim |x|^{-n}
.\end{multline*}
Next,
$$
\int_{\delta|x|}^{|x|} \fra\Theta_r(x)(- q'(r))\,dr
\le \int_{\delta|x|}^{|x|} (- q'(r))\,dr\le q(|x|)+q(\delta|x|)
.$$
Thus the asymptotics of our integral is $\sim |x|^{-n}$.
\hfill$\square$

\sm

We also need estimates of 
the partial derivatives of functions $\fra\phi(x)$.
Let $e_j$ be a basis vector in $\R^{2n}$
$$
\frac{\partial}{\partial x_j}
\fra \phi_j
=\frac d{ds}\int_{\cAL_{x+se_j}}\phi(L)\,d\sigma_{x+se_j}(L) \biggr|_{s=0}=
\frac d{ds}\int_{\cAL_{x}}\phi(-se_j+L)\,d\sigma_{x}(L)\biggr|_{s=0}
.$$ 
We have 
$$\frac d{ds}\phi(-se_j+L)\in\cS(\cAL(n,\R))\Bigr|_{s=0}\in\cS(\cAL(n,\R)).
$$
by definition of $\cS(\cAL(n,\R))$,
and we get the same asymptotic.

This finishes the proof of Theorem \ref{th:distributions}.

\section{The image of $\cR_{flat}$}

\subsection*{Fourier transform}
Denote by $\Symm_0(n+1)$ the space  of real symmetric  $(n+1)\times(n+1)$ matrices 
of the form
$\begin{pmatrix}A&\alpha
\\
\alpha^t &0
\end{pmatrix}$. 
 The group $\GL(n,\R)$ acts on $\Symm_0(n+1)$ by transformations
 $$
 \begin{pmatrix}A&\alpha
\\
\alpha^t &0
\end{pmatrix}\mapsto  
\begin{pmatrix}g&0\\0&1 \end{pmatrix}
\begin{pmatrix}A&\alpha
\\
\alpha^t &0
\end{pmatrix}
\begin{pmatrix}g&0\\0&1 \end{pmatrix}^t
. $$

It will be convenient for us to present coordinates $(T,\tau)$ on $\cAL(n,\R)$ in the matrix form
$$\cT=\begin{pmatrix}
  T&\tau\\ \tau^t&0
 \end{pmatrix}\in \Symm_0(n+1).$$

We write the Fourier transform on $\Symm_0(n+1)$ by
$$
G(\cX)=\int_{\Symm_0(n+1)} F(\cT)e^{\frac i2 \,\tr \cT\cX}\,d\cT,\qquad
\text{where $\cX=\begin{pmatrix} X&\xi\\\xi^t&0\end{pmatrix}\in \Symm_0(n+1)$.} 
$$

\begin{lemma}
Let $F$ be a tempered distribution  on $\Symm_0(n+1)$ satisfying the determinantal system of equations.
Then its Fourier transform $\phi$ is supported by  
matrices of rank $\le 2$.
\end{lemma}

{\sc Proof.}
Transforming the determinantal system  to the language of the Fourier transform, we get
\begin{equation}
P\cdot \phi=0
,
\label{eq:det-var}
\end{equation}
where $P$ ranges in $3\times 3$ minors of the matrix
\begin{equation}
\cX=
\begin{pmatrix}
x_{11}&\dots& x_{1n}&\xi_1
\\
\vdots&\ddots&\vdots &\vdots
\\
x_{1n}&\dots &x_{nn}&\xi_n
\\
\xi_1&\dots&\xi_n&0 .
\end{pmatrix}
\label{eq:det-var2}
\end{equation}

\subsection*{Matrices of rank 2.}
Consider  the set $\Delta\subset \Symm_0(n+1)$ consisting of all  matrices
of rank $\le 2$.

\begin{proposition}
{\rm a)} The set $\Delta$ is a union of two closed sets:
\begin{itemize}
	\item 
The $(2n-1)$-dimensional  set $M$, consisting of all matrices of the form
 $$\begin{pmatrix}
           X&0\\0&0
          \end{pmatrix},\qquad \rk X\le 2
          ;$$
\item
The $2n$-dimensional set $N$ consisting of all matrices of the form
\begin{equation}
\begin{pmatrix}
 y^tz+z^t y&z
\\
z^t&0
\end{pmatrix},
\label{eq:yz}
\end{equation}
where $z$,
$y$ - are vector-columns.
\end{itemize}

The intersection $M\cap N$ consists of all matrices $\begin{pmatrix}
           X&0\\0&0
          \end{pmatrix}$
          such that $X$ is sign-indefinite of rank $2$ or $\rk X\le 1$.
          
          \sm
          
{\rm b)}   There is a unique {\rm(}up to a factor{\rm)}  $\SL(n,\R)$-invariant measure on $N$, it is given by
$$
dy_1\dots dy_n \, dz_1\dots dz_n
.$$
\end{proposition}

{\sc Remark.}
An elementary model of such situation is the Whitney umbrella,
i.e., the surface $z(x^2+y^2)=xy$ in $\R^3$. It consists
of the 2-dimensional ruled surface $z=\frac{x^2-y^2}{x^2+y^2}$ in the layer $-1\le z\le 1$ and the 
line $x=y=0$. Only the segment $-1\le z\le 1$ of the line $x=y=0$ is contained in the closure 
of the surface (this segment is a line of self-intersection).
\hfill $\square$

\sm

{\sc Proof.} a)
We reduce $X$ to diagonal form. The only possible variants of diagonal values are 
$(\pm1,\pm 1,0,0\dots)$, $(\pm1,0,0\dots)$, $(0,0,0,\dots)$.
Next,
$$
\det\begin{pmatrix} 1&0&\xi_1\\
0&0&\xi_2\\
\xi_1&\xi_2&0
\end{pmatrix}=\xi_2^2,
\quad
\det\begin{pmatrix} 1&0&\xi_1\\
0&\pm 1&\xi_2\\
\xi_1&\xi_2&0
\end{pmatrix}=\xi_2^2\pm\xi_2^2
$$
Let diagonal values of $X$ be  $\pm(1,1,0,0,\dots)$.
Then we see that  $\xi=0$.

Next, assume that $\xi\ne0$. By applying an appropriate element
of $\GL(n,\R)$ we reduce it to the form
$\begin{pmatrix}
 1\\0\\\vdots 
 \end{pmatrix}
$.
The only way to extend $\xi$  to a matrix $\cX\in\Symm_0(n+1)$ of rank $\le 2$
is
$$
\cX=
\begin{pmatrix}
 2b_1&b_2&b_3&\dots& 1
\\
b_2&0&0&\dots& 0
\\
b_3&0&0&\dots& 0
\\
\vdots&\vdots& \vdots&\ddots &\vdots
\\
1&0&0&\dots& 0
\end{pmatrix}
.$$
We get an expression of form (\ref{eq:yz}) with
$$z=
\begin{pmatrix}
 1\\0\\\vdots\\0 
 \end{pmatrix},\qquad
 y=\begin{pmatrix}
 b_1\\b_2\\\vdots\\b_n 
 \end{pmatrix}
 .\qquad\square
$$

%
%


\sm

b) The group $\SL(n,\R)$ acts on the manifold $N$
as
$ y\mapsto gy, z\mapsto gz$.
\hfill $\square$

\subsection*{Description of the image of the Radon transform.}
Consider a function $\chi(u,v)\in\cS(\R^{n}\times \R^{n})$.
For a function $\phi\in\cS(\Symm_0(n+1))$ we consider the value
\begin{equation}
\phi\mapsto
\int_{\R^{n}\times\R^{n}} \phi
  \begin{pmatrix}
 y^tz+z^t y&z
\\
z^t&0
\end{pmatrix}
                                \chi(u,v)\,du\,dv
\label{eq:distribution-rank2}
.
\end{equation}

In fact, this is  a tempered distribution  supported by $N$.

\begin{theorem}
The following conditions are equivalent:
\begin{itemize}
	\item 
$F$ is contained in the image of $\cS(\R^{2n})$ under the Radon transform.
\item
$F$ is a Fourier transform of a distribution {\rm(\ref{eq:distribution-rank2})}.
\end{itemize}
\end{theorem}

{\sc Proof.}
Consider the Fourier transform of our distribution,
\begin{multline*}
G
\begin{pmatrix}
T&\tau\\ \tau^t&0
\end{pmatrix}
=\int_{\R^n\times \R^n}
\chi(x,z)e^{\frac i2 \tr T(xz^t+xy^t )+i\,z^t\tau}\,dx\,dz
=\\=
\int_{\R^n\times \R^n} \chi(x,z)e^{i(z^t Tx+z^t \tau) }\,dx\,dz
.
\end{multline*}

This expression is the Radon transform of the function
\begin{equation}
\psi(\xi,z)=\int_{\R^n} \chi(x,z) e^{i\xi^t x}\,dx
\label{eq:partial-fourier}
.\end{equation}
Note that the  operator (\ref{eq:partial-fourier})
is a bijection $\cS(\R^n\times\R^n)\to \cS(\R^n\times\R^n)$.

\subsection*{Extraneous solution}
Let $G$ be a tempered distribution on $\Symm_0(n+1)$ supported by $\Delta$
and satisfying the equations
(\ref{eq:det-var})--(\ref{eq:det-var2}). Notice that any measure supported by $\Delta$ satisfies these
equations.  On the other hand, the algebraic variety $\Delta$
is singular, therefore  some partial derivatives of measures
also are admissible.

Then the Fourier transform of $G$ is a solution of the determinantal system on $\Symm_0(n+1)$.
However, subsets  $M$, $N\subset \Delta$ produce two different families of such solutions.
We also observe that the image of the Radon transform is not dense
in the topology of $\cS'(\Symm_0(n+1))$ in the space of solutions of the determinantal system.

\section{Decomposition in tensor product}

\subsection*{Weil representation}
Recall that the {\it Weil representation} $\on{We}$ of $\Sp(2n,\R)$ is given as follows
(see, e.g., \cite{Neretin96} or \cite{Neretin11}, Section 2.1).
For generators of $\ASp(2n,\R)$ we have
\begin{align}
\on{We}
\begin{pmatrix} (a^t)^{-1} & 0&0\\ 0 & a&0\\0&0&1\end{pmatrix} 
f(u)&
=
 \pm \det(a)^{1/2} f(a u )
\label{eq:we1}
\\
\on{We}
\begin{pmatrix}1 & b&0\\ 0 & 1&0\\0&0&1\end{pmatrix}
f(u)&=f(u)e^{\frac i2  u^t b u}
\label{eq:we2}
\\
\on{We}
\begin{pmatrix}
0&1&0\\-1&0&0\\0&0&1
\end{pmatrix} f(u)&=
 \left( \frac{i}{2\pi}\right)^{n/2} \int e^{i\xi^t u}f(\xi)\,d\xi
 \label{eq:we3}
 \\
 \on{We}
\begin{pmatrix}1 & 0&r\\ 0 & 1&0\\0&0&1\end{pmatrix}
f(u)&= f(u+r)
\label{eq:we4}
\\
\on{We}
\begin{pmatrix}1 & 0&0\\ 0 & 0&s\\0&0&1\end{pmatrix}f(u)
&
= f(u)e^{i s^t b }
\label{eq:we5}
\end{align}
The group of operators generated by (\ref{eq:we1})--(\ref{eq:we5})
is isomorphic to a central extension of
$\ASp(2n,\R)$ by $\T=\R/2\pi Z$.
The central subgroup $\T$ acts by multiplications by
$e^{i\theta}$, $\theta\in\R$.

In fact we get a
projective representation
of $\ASp(2n,\R)$, 
\begin{equation}
\on{We}(g_1)\on{We}(g_2)=\sigma(g_1,g_2)\cdot \on{We}(g_1,g_2)
,
\label{eq:projective}
\end{equation}
where $\sigma(g_1,g_2)$ is a scalar factor, $|\sigma|=1$.
There are two reasons for the appearance of this factor:
\begin{itemize}
	\item First, there is the sign $\pm$ in (\ref{eq:we1}). In fact, the operators 
        (\ref{eq:we1})--(\ref{eq:we3}) determine a representation of a double covering of 
        $\Sp(2n,\R)$.
  \item Secondly, the operators (\ref{eq:we3}) and (\ref{eq:we4}) do not commute, they generate the Heisenberg group.
\end{itemize}

\newpage
\subsection*{Tensor product}

\begin{theorem}
The representation of $\ASp(2n,\mathbb R)$ on $\mathcal S(\mathbb R^{2n})$ and on 
$L^2(\mathbb R^{2n})$ is isomorphic to the tensor product 
of the Weil representation $\on{We}$ with its 
complex conjugate representation $\ov {\on{We}}$.
The intertwining unitary  operator $L^2(\R^{2n})\to L^2(\R^n)\otimes L^2(\R^n)$ is given by
$$
H f(u,v)=\frac{1}{(2\pi)^n}\int_{\R^n} f\Bigl(x,\frac {u-v}{\sqrt 2}\Bigr) \exp\Bigl\{\frac{i}{\sqrt 2}  x^t(u+v)\Bigr\} \,dx
$$
This operator also is a bijection $\cS(\R^{2n})\to \cS(\R^n)\hat\otimes \cS(\R^n)$.
\end{theorem}

{\sc Remarks.} a)
We recall that the  tensor product  $L^2(M)\otimes L^2(M)$
in the category of Hilbert spaces
is canonically isomorphic to
$L^2(M\times N)$.

\sm

b) 
The operator $H$ is a composition of the Fourier transform
with respect to half of variables, 
$$
\cF f(x,y)=\frac{1}{(2\pi)^n}\int_{\R^n} f(x,y)\,e^{i  x^t \xi}\,dx
$$
and of the rotation of the plane $(y,\xi)$ by $45$ degrees. 

\sm

c) We have
$$
\bigl(\on{We}(g_1)\otimes \ov{\on{We}}(g_1)\bigr) \cdot \bigl(\on{We}(g_2)\otimes \ov{\on{We}}(g_2)\bigr)
=\sigma(g_1,g_2)\ov{\sigma(g_1,g_2)}\bigl(\on{We}(g_1g_2)\otimes \ov{\on{We}}(g_1g_2)\bigr)
.$$
Since $|\sigma|=1$, the scalar factor disappears and we get a linear representation.
\hfill $\square$

{\sc Proof} is straightforward, we  look at the images of the operators
(\ref{eq:we1})--(\ref{eq:we5}) under the partial Fourier transform
and the rotation.
\hfill $\square$

\section{Pair of wedges}

\subsection*{Siegel half-plane}
Denote by $\cZ^+(n)$ the space of $n\times n$ matrices $P$ 
such that its imaginary part $\Im P$ is positive-definite (this set is called {\it Siegel upper half-plane}).
By $\cZ^-(n)$ we denote the space of matrices with negative-definite $\Im P$.
 The group
$\Sp(2n,\R)$ acts on $\cZ^+(n)$ by linear fractional transformations
\begin{equation}
P\mapsto (aP+b)(cP+d)^{-1}
\label{eq:siegel-action}
\end{equation}
Note that a complex matrix $P$ determines a point of $\cL(n,\C)$ and  formula
(\ref{eq:siegel-action}) is a special case of (\ref{eq:perenos-2}). Note that this expression is
 defined everywhere.

More generally, consider the set  
$\cA\cZ^+(n):=\cZ^+(n)\times \C^n$, denote its points by $(P,\pi)$.
 The group
$\ASp(2n,\R)$ acts on this space by transformations
\begin{equation}
(P,\pi)\mapsto
\\ \bigl((aP+b)(cP+d)^{-1},\quad  a\pi+r\, -(aP+b)(cP+d)^{-1} (c\pi+s)  \bigr)
\label{eq:again}
\end{equation}

\subsection*{Geometric interpretation}
 As we have seen in Section 2 (see (\ref{eq:form-Lambda}), (\ref{eq:M})) the complex space
 $\C^{2n}$ is equipped with two forms,
 the skew-symmetric form $\Lambda$ and the indefinite
 Hermitian form $M$. Consider the action of the real group
 $\Sp(2n,\R)$ on the complex Lagrangian Grassmannian
 $\cL(n,\C)$. Evidently, the inertia indices of $M$ on $L\in\cL(n,\C)$ are invariants 
 of this action. Denote by $\cL(n,\C)_+$ (resp $\cL(n,\C)_-$)  the set of Lagrangian subspaces strictly positive 
 (resp. negative)
 with respect to $M$. Denote by $\cAL(n,\C)_\pm$ 
 the set of subspaces of form $\xi+L$, where $\xi\in\C^{2n}$, $L\in\cL(n,\C)_\pm$.
 
 We represent such subspace as a graph of a map $W_-\to W_+$
 as in (\ref{eq:aff-T-tau}) and get  formula (\ref{eq:again}).

\subsection*{Another realization of the Weil representation}
Let $f$ be a function on $\R^n$. We define an integral transform
$\frK$ by 
\begin{equation}
\frK f(P,\pi):=
\frac{1}{(2\pi)^{-n/4}}\int_{\R^n} f(u)e^{\frac i2 u^t Pu+\frac i{\sqrt{2}}u^t\pi}\,du
.
\label{eq:xPx}
\end{equation}
Obviously, for any $f\in \cS'(\R^n)$ we get a holomorphic function
on $\cA\cZ^+(n)$.

This operator sends the space $L^2(\R^{n})$ to a certain Hilbert $\cH$ space of holomorphic
functions on $\cA\cZ^+(n)$. We wish to discuss this space and the action of 
 of the group $\ASp(2n,\R)$ in this space.

\subsection*{Reproducing kernel.} First, our Hilbert space
can be  described by the standard language of reproducing kernels.
Consider a Hilbert space $H$ whose elements are continuous functions $f$
on a domain $\Omega\subset\R^N$. Assume that the point evaluations $f\mapsto f(a)$
are continuous linear functionals on $H$ for all $a\in\Omega$. Then for each
$a\in\Omega$ we have
$$
f(a)=\la f,\phi_a\ra_H
$$
for some $\phi_a\in H$.
The function
$$
K(a,b)=\la\phi_a,\phi_b\ra=\phi_a(b)=\ov{\phi_b(a)}
$$
is called {\it the reproducing kernel} of the space $H$. 
The space $H$ (the set of function and the inner product) is uniquely determined by the reproducing
kernel; see, e.g., \cite{Neretin11}, Section 7.1.

\begin{proposition}
In our case,
\begin{multline}
K(P,\pi;\,R,\rho)=\\=
\phi_{R,\rho}(P,\pi)
\det[(-i(P-\ov R))^{-1/2}]
\exp\left\{-i (\pi-\ov \rho)^t(P-\ov R)^{-1} (\pi-\ov \rho) \right\}
,
\label{eq:K=repr}
\end{multline}
where we choose the holomorphic branch of the determinant satisfying
$$
\det[(-i(P-\ov R))^{-1/2}]\Bigr|_{P-\ov Q=i}=\det[(-i\cdot i)^{-1/2}]=1
.
$$

\end{proposition}

{\sc Proof.} We 
set
$$
\Phi_{R,\rho}(x)=\frac 1{2\pi^{n/2}}
\exp\bigl\{\frac i2 x^t \ov Rx+\ov\rho^t x\bigr\}
$$
It is contained 
in $L^2(\R^n)$ if $\Im P>0$.
Then
$$
\frK\Phi_{R,\rho}(P,\pi) =\phi_{R,\rho} (P,\pi)  
$$
$$
\frK f(R,\rho)=\la f, \Phi_{R,\rho}\ra_{L^2(\R^n)}=
\la \frK f, \phi_{R,\rho}\ra_{\cH}
.
$$
This implies the reproducing kernel property.
\hfill $\square$

\subsection*{Differential equations}

\begin{theorem}
The functions  $\frK f(P,\pi)\in \cH(\cA\cZ^+(n))$
satisfy the equations
$$
D \,\frK \,F(P,\pi)=0
,$$
where $D$ ranges in all the $2\times 2$ minors of the matrix
$$
\begin{pmatrix}
2\frac\partial{\partial p_{11}} & \frac\partial{\partial p_{12}} & \dots&
     \frac\partial{\partial p_{1n}} &\frac\partial{\partial \pi_1}
\\
\frac\partial{\partial p_{12}} & 2\frac\partial{\partial p_{22}} & \dots&
   \frac\partial{\partial p_{2n}} &\frac\partial{\partial \pi_2}
\\
\vdots&\vdots& \ddots&\vdots&\vdots
\\
\frac\partial{\partial p_{1n}} & \frac\partial{\partial p_{2n}} & \dots &
          2 \frac\partial{\partial p_{nn}} &\frac\partial{\partial \pi_n}
\\
\frac\partial{\partial \pi_1}&\frac\partial{\partial \pi_2}&\dots& 
\frac\partial{\partial \pi_n}&i/2
\end{pmatrix}
\label{eq:system-Z}
$$
\end{theorem}

{\sc Proof.} a) Indeed the functions $\phi_x(P,\pi)=e^{\frac i2x^t Px+\frac i{\sqrt 2}\pi^t x}$
satisfy these equations. Next, we differentiate the integral (\ref{eq:xPx})
 with respect to the parameters $P$, $\pi$ and get the required statement.
However, some justification here is necessary. For that we refer
to Theorem 7.7.6 of \cite{Neretin11}. \hfill$\square$

\sm

Note that the image of any tempered distribution under the transform $\frK$
satisfies the same system of differential equations.
 
 \subsection*{The action of the group $\ASp(2n,\R)$ in the space of holomorphic functions.} 
First, we transfer the operators of
$\Sp(2n,\R)$

\begin{lemma}
 The group $\ASp(2n,\R)$ acts in the space $\cH$ by the formulas
\begin{align}
S
\begin{pmatrix} (a^t)^{-1} & 0&0\\ 0 & a&0\\0&0&1\end{pmatrix} 
F(P,\pi)&=\pm\det(a)^{-1/2}(a)F((a^t)^{-1}P a^{-1},(a^t)^{-1}\pi)
\label{al:1}
\\
S\begin{pmatrix} 1 & b&0 \\ 0 & 1&0\\0&0&1\end{pmatrix} 
F(P,\pi)&=F(P+b,\pi)
\label{al:2}
\\
S
\begin{pmatrix} 0 & 1 &0\\ -1 & 0&0\\0&0&1\end{pmatrix} F(P,\pi)
&=\det(-P)^{-1/2}e^{-\frac i4 \pi^tP^{-1}\pi} F(-P^{-1}, P^{-1}\pi)
\label{al:3}
\\
S\begin{pmatrix} 1 & 0&0 \\ 0 & 1&0\\0&s&1\end{pmatrix}  F(P,\pi)&= e^{\frac i4 s^tPs -\frac 1{\sqrt 2}s^t\pi} F(P,\pi-Ps)
\label{al:4}
\\
S\begin{pmatrix} 1 & 0&0 \\ 0 & 1&0\\r&0&1\end{pmatrix}  F(P,\pi)&=F(P,\pi+r)
\label{al:5}
\end{align}
\end{lemma}

This is proved by straightforward calculations.

\subsection*{Product of wedges.}
We can apply the same construction to the complex conjugate  of 
the Weil representation.
We consider the integral transform
$$
\frK f(P,\pi):=\frac{1}{(2\pi)^{-n/4}}\int_{\R^n} f(u)e^{-\frac {i}2 u^t Pu-
\frac i{\sqrt{2}}u^t\pi}\,du,
$$
an operator from $cS'(\mathbb R^n)$
to the space $\cH_-$ of holomorphic functions on
$\cA\cZ^-(n)$. All formulas are the same, only $i$ must be replaced by $-i$ everywhere. 

Our representation of $\ASp(2n,\R)$ is a tensor product of the Weil representation and its complex
conjugate representation. Therefore we can realize our representation
 in the space $\cH^+\otimes \cH^-$ of
 holomorphic functions 
on $\cA\cZ^+(n)\times \cA\cZ^-(n)$.

\subsection*{Intertwining operator.} 

\begin{theorem}
\label{th:intertwining}
The operator intertwining the representation of $\ASp(2n,\R)$ in $L^2(\R^{2n})$
and in $\cH^+\otimes \cH^-$ is given by the formula
{\small
\begin{multline*}
\frJ f(P,\pi;Q,\kappa)
= 
\frac 1{(2\pi)^n}
\det[(-i(Q-P))^{-1/2}]
\int_{\R^{2n}} f(x,y) \cdot
\\
\exp\left\{
-\frac i2 \begin{pmatrix}
(x+\pi+\kappa)^t& y^t
\end{pmatrix}
\begin{pmatrix}
2(P-Q)^{-1}& (P-Q)^{-1}(P+Q)
\\
(P+Q)(P-Q)^{-1}& -2 (P^{-1}-Q^{-1})
\end{pmatrix}
\begin{pmatrix}
x+\pi+\kappa\\ y
\end{pmatrix}
\right\}
\\
\cdot\,dx\,dy
\end{multline*}
}
 \end{theorem}

We omit a formal proof. It is sufficient to verify that the operator commutes
with the generators of $\ASp(2n,\R)$.

\subsection*{Additional remarks}
We note that the matrix
$$g= i\begin{pmatrix}
2(P-Q)^{-1}& (P-Q)^{-1}(P+Q)
\\
(P+Q)(P-Q)^{-1}& -2 (P^{-1}-Q^{-1})
\end{pmatrix}$$
is complex symplectic. Moreover, the
eigenvalues of $\begin{pmatrix} 0&1\\-1&0\end{pmatrix} g$
 are $\pm 1$, and 
the eigenspaces
are elements of $\cAL(n,\C)_+$ and $\cAL(n,\C)_-$ with coordinates $P$, $Q$.

The Radon transform corresponds to the limit case $Q=P$
(this case corresponds to the common edge of two wedges).

\end{document}